\documentclass[journal]{IEEEtran}

\usepackage{cite}
\usepackage{color}

\ifCLASSINFOpdf
   \usepackage[pdftex]{graphicx}
\else
\fi

\usepackage{amsmath}
\usepackage[ruled,vlined]{algorithm2e}

\begin{document}

\title{Optimal Distributed Energy Resource Coordination: {\color{black} A} Decomposition {\color{black} Method} Based on Distribution Locational Marginal Costs}

\author{Panagiotis~Andrianesis,~\IEEEmembership{Member,~IEEE},
        Michael~Caramanis,~\IEEEmembership{Senior~Member,~IEEE}
        and~Na~Li,~\IEEEmembership{Member,~IEEE}
\thanks{P. Andrianesis and M. Caramanis are with the Systems Eng. Div., Boston University, Boston, MA: panosa@bu.edu, mcaraman@bu.edu.
N. Li is with the School of Engineering and Applied Sciences, Harvard University, Cambridge, MA 02138 USA: nali@seas.harvard.edu.
Research partially supported by Sloan Foundation G-2017-9723, NSF AitF 1733827, and PGMO-IROE 2020-0009.
}
}

\maketitle
\begin{abstract}
In this paper, 
we consider the day-ahead operational planning problem of a radial distribution network 
hosting Distributed Energy Resources (DERs) 
including rooftop solar and storage-like loads, 
such as electric vehicles. 
We present a novel 
decomposition {\color{black} method} 
that is based on a centralized AC Optimal Power Flow (AC OPF) problem
interacting iteratively with self-dispatching DER problems
adapting to real and reactive power Distribution Locational Marginal Costs.
We illustrate the applicability and tractability of the proposed {\color{black} method} 
on an actual distribution feeder, 
while modeling the full complexity of spatiotemporal DER capabilities and preferences,
and accounting for instances of non-exact AC OPF convex relaxations.
We show that the proposed {\color{black} method} achieves optimal Grid-DER coordination,
by successively improving feasible AC OPF solutions,
and discovers spatiotemporally varying marginal costs in distribution networks
that are key to optimal DER scheduling by modeling losses,
ampacity and voltage congestion,
and, most importantly, dynamic asset degradation.
\end{abstract}


\IEEEpeerreviewmaketitle

\section{Introduction}

\IEEEPARstart{I}{ntegration} of increasingly prevalent Distributed Energy Resources (DERs), 
including Solar Photovoltaic (PV) and storage-like loads, such as Electric Vehicles (EVs), 
presents major challenges with notable impacts on local costs and constraints, 
e.g., losses, under/over voltage, infrastructure congestion, transformer degradation 
\cite{EVproject2013, SMUD}. 
However, 
it is also associated with positive synergies 
that may provide significant benefits
by putting the DER mostly underutilized inverter capacity to dual use
for the provision of services to the Distribution Grid,
e.g., Volt/VAR control. 
The day-ahead operational planning problem 
of a distribution network hosting DERs 
is key to achieving optimal Grid-DER coordination.  

Unlike current operational planning transmission system practice that relies on 
linear DC Optimal Power Flow (OPF) models, 
in distribution networks, accurate modeling of real/reactive power flows, voltages and currents is
more important, thus requiring AC OPF formulations.  
Earlier literature exploring distribution network markets 
has adhered to the transmission paradigm relying on DC OPF and
uniform price-quantity bids \cite{LiEtAl_2014, HuangEtAl_2015, BaiEtAl-DLMP}.
We argue that, although these approaches may have derived useful results,  
they have not gone far enough towards the discovery of system-wide spatiotemporal marginal costs 
that are key to the efficient dispatch of DERs accounting for their preferences and degrees of freedom \cite{CaramanisEtAl2016}.

Given the almost ubiquitous radial topology used in distribution networks,
the ``branch-flow'' model (\emph{a.k.a.} DistFlow), 
introduced by \cite{BaranWu1989},
has been recently revisited \cite{FarivarLow2013}
and employed widely. 
Using a relaxed branch flow model, 
the AC OPF problem can be formulated as a Second Order Cone Programming (SOCP) problem. 
Current exploration of the SOCP formulation \cite{KocukEtAl16} has enabled recent works 
to investigate convex relaxations and the conditions 
under which they are exact 
\cite{GanEtAl15, Abdelouadoud_EtAl_2015, HuangEtAl_2017, WeiEtAl_2017, NickEtAl_2018, YuanEtAl-DLMP}. 
Notably, however, 
real systems experiencing increasing distributed generation
exhibit flow reversals 
that violate the conditions guaranteeing exact convex relaxations 
and lead to erroneous solutions that violate the physics of power flow. 
Furthermore, negative LMPs, 
neither commonplace but at the same time not unlikely either,
create conditions under which these relaxations are not exact. 

We note that DERs, despite their huge numbers and complex preferences including spatiotemporal dynamics, 
are amenable to modern decomposition approaches involving parallelizable distributed algorithms \cite{MolzahnEtAl_2017}, 
such as the alternating direction method of multipliers (ADMM) described in \cite{Kraning14}. 
Although decomposition approaches date several decades back, 
e.g., ADMM was first proposed in the 70’s (a brief historical review is provided in \cite{Ber}), 
they have recently benefited from renewed attention, among others, in the work of \cite{BoydEtAl10}.
Despite the significant attention that distributed solution approaches to the AC OPF problem
have attracted, 
their application to real systems has remained a challenge, since 
they have typically assumed convexity, and, most importantly, 
they do not provide a physically meaningful solution until they converge.
Industry has therefore been reluctant to adopt them.

In our recent work \cite{IEEEtsg20}, 
we presented a novel centralized AC OPF model 
encompassing transformer degradation as a short-run network variable cost, 
and moreover additive real/reactive power Distribution Locational Marginal Cost (DLMC) components related to: 
\emph{(i)} the costs of real/reactive power transactions at the T\&D interface; 
\emph{(ii)} real/reactive power marginal losses; 
\emph{(iii)} voltage/ampacity congestion, and 
\emph{(iv)} a new transformer degradation marginal cost component. 
Our detailed transformer degradation model captured the impact 
of incremental transformer loading on its top oil temperature during a specific time period, 
impacting not only its Loss of Life (LoL) during that period, 
but also during subsequent time periods.
We further employed real distribution feeders 
to exemplify the use of DLMCs as financial incentives 
that convey sufficient information to optimize 
both Distribution Network and DER (PV, EV) operation, 
which also extends to a yearly planning horizon. 
As such, 
we considered DERs as non-wires alternatives 
to the distribution network investment planning problem \cite{hiccs, LMV}.

In this paper, we present a novel 
decomposition {\color{black} method} 
that utilizes a centralized AC OPF network-optimization problem 
interacting with multiple DER-specific self-dispatch problems 
that adapt DER schedules to real and reactive power DLMCs generated by the network-optimization problem.
This adds to the literature \cite{MolzahnEtAl_2017} 
a novel algorithm enabling the derivation of 
mutually adapted spatiotemporal marginal costs 
and optimal DER schedules. 
Our approach promises to create a new learning curve for DER owners and distribution utilities. 
Striking a balance between centralized control and distributed self-dispatch, 
the proposed 
decomposition is scalable 
and enables massive DER participation to the grid 
with high modeling fidelity of their preferences, dynamics, and degrees of freedom,
which, among others, improves the accuracy of demand forecasts available to transmission level Wholesale Markets.
The applicability and tractability of the proposed {\color{black} method} is further explored
in actual distribution feeder case studies
where the full complexity of spatiotemporal DER capabilities and preferences is modeled.

Preliminary results have been presented in \cite{Allerton},
where we focused on interpreting DLMCs as price signals, 
provided an informal sketch of the idea,
and illustrated it on two service transformers described in \cite{IEEEtsg20}.
This paper extends our contribution as follows:
First, 
we provide a concise 
decomposition-based modeling framework
with a stylized proof/description of the underlying proximal algorithm,
{\color{black}and highlight the key differences of our proposal with popular decomposition algorithms, 
such as dual decomposition and ADMM.}
Second, 
we complement the framework
by embedding the idea of \cite{WeiEtAl_2017} in the proposed {\color{black} method}
considering cases where convex relaxations are not exact.
Third, 
we illustrate the applicability of our {\color{black} method} on an actual distribution feeder, 
validate the proof-of-concept with numerical results, 
{\color{black}and discuss practical considerations.}

The remainder of the paper is organized as follows.
Section \ref{OPF} provides the centralized formulation of the enhanced AC OPF problem.
Section \ref{Coordination} presents the main idea of the proposed {\color{black} method} 
for optimal Grid-DER coordination. 
It describes the 
decomposition, 
its application to the enhanced AC OPF problem, 
and the proposed remedy for problem instances with non-exact SOCP relaxations.
Section \ref{Numerical} presents numerical results 
implementing the proposed method 
on a case study involving an actual distribution feeder.
{\color{black} Section \ref{Discuss} discusses practical considerations and contributions of the method.} 
Section \ref{Conclusions} concludes and provides further research directions.

\section{Enhanced AC OPF Problem} \label{OPF}

In this section, we briefly present the AC OPF operational planning problem 
in a radial distribution network hosting EVs and PVs, 
enhanced to include transformer degradation \cite{IEEEtsg20}.

Let $\mathcal{N}= \{ 0,1,...,N \}$ denote the set of nodes, 
where node $0$ represents the root node.
Let $u(j)$ be the unique upstream/parent node of node $j$,
and $\mathcal{C}_j$ the set of the children of node $j$.
The root node $0$ has no parent, 
whereas leaf nodes have no children.
The line connecting nodes $u(j)$ and $j$ is also denoted by $j$,
so that $\mathcal{N^+} \equiv \mathcal{N} \setminus \{0\}$ also represents the set of lines.
Each line $j$ is characterized by resistance $r_{j}$ and reactance $x_{j}$.
Let $T$ be the length of the optimization horizon,
and time period $t \in \mathcal{T}^+$, with $\mathcal{T} = \{0,1,...,T \}$, $\mathcal{T^+} \equiv \mathcal{T} \setminus \{ 0 \}$. 
For line $j$, and time period $t$,
$l_{j,t}$ denotes the magnitude squared current, 
and $P_{j,t}$ and $Q_{j,t}$ the sending-end real and reactive power flow, respectively.
Assuming, 
without loss of generality, 
that the root node has only one child (node 1), 
the net injections at the root node should equal the real and reactive power flows $P_{1,t}$ and $Q_{1,t}$.
For node $j$, and time period $t$,
$v_{j,t}$ denotes the magnitude squared voltage,
and $p^d_{j,t}$ and $q^d_{j,t}$ the conventional demand for real and reactive power, respectively.
Service transformers are a subset of lines, denoted by $\mathcal{N}^+_D$,
with the downstream node being a leaf node.
For service transformer $j \in \mathcal{N^+_D}$, and time period $t$,
$D_{j,t}$ denotes the degradation (LoL) per time period,
and $h_{j,t}$ the top-oil temperature.
EVs and PVs are denoted by sets $\mathcal{E}$, and $\mathcal{S}$, respectively. Subset $\mathcal{E}_{j,t} \subset \mathcal{E}$ denotes EVs connected at node $j$, during time period $t$,
and $\mathcal{S}_j \subset \mathcal{S}$ denotes PVs connected at node $j$.
Let $p_{e,t}$ and $q_{e,t}$ denote the real and reactive power consumption of EV $e \in \mathcal{E}$.
Let $p_{s,t}$ and $q_{s,t}$ denote the real and reactive power provision of PV $s \in \mathcal{S}$.

The enhanced AC OPF problem aims at selecting the DER schedules and network variables
to minimize the aggregate real and reactive power cost and the transformer degradation cost, 
where $c_t^P$ is typically the LMP at the T\&D interface (substation), 
$c_t^Q$ is the opportunity cost for the provision of reactive power, 
$c^D_j$ is an hourly transformer cost representing the cost of losing one hour of transformer life.
Unless otherwise mentioned, $j \in \mathcal{N}^+, t \in \mathcal{T}^+$.
\begin{equation} \label{Obj2}
\text{min}
\sum_{t} c^P_t P_{1,t}
+ \sum_{t} c^Q_t Q_{1,t} 
+ \sum_{j \in \mathcal{N}_L,t} c^D_j D_{j,t} ,
\end{equation}
\emph{subject to}: 

\underline{Power Flow Constraints}: 
\begin{align} \label{EqRealBalance}
\begin{split}
P_{j,t} - r_{j} l_{j,t} & = \sum_{j' \in \mathcal{C}_j} P_{j',t} + p^d_{j,t} \\ 
& + \sum_{e \in \mathcal{E}_{j,t}}  p_{e,t} - \sum_{s \in \mathcal{S}_j} p_{s,t} \rightarrow (\lambda_{j,t}^P), \quad \forall j,t,
\end{split}
\end{align}
\begin{align} \label{EqReactiveBalance}
\begin{split}
Q_{j,t} - x_{j} l_{j,t} & = \sum_{j' \in \mathcal{C}_j} Q_{j',t} + q^d_{j,t} \\ 
&+ \sum_{e \in \mathcal{E}_{j,t}}  q_{e,t} - \sum_{s \in \mathcal{S}_j} q_{s,t} \rightarrow (\lambda_{j,t}^Q), \quad \forall j,t,
\end{split}
\end{align}
\begin{equation} \label{EqVoltageDef}
v_{j,t} = v_{u(j),t} - 2 r_{j} P_{j,t} - 2 x_{j} Q_{j,t} + \left( r_{j}^2 + x_{j}^2 \right) l_{j,t}, \,\,\, \forall j,t,
\end{equation}
\begin{equation} \label{EqCurrentDef}
 v_{u(j),t} l_{j,t} \geq P_{j,t}^2 + Q_{j,t}^2, \quad \forall j,t.
\end{equation}
Constraints \eqref{EqRealBalance} and \eqref{EqReactiveBalance} define the real and reactive power balance, respectively.
Constraint \eqref{EqVoltageDef} defines the voltage drop along a line.
Inequality \eqref{EqCurrentDef} is the SOCP relaxation 
--- introduced by \cite{FarivarLow2013} for the DistFlow model in \cite{BaranWu1989} --- 
of the (non-convex) equality constraint that defines apparent power.
\begin{align}
& \text{\underline{Voltage Limits}: } & \quad \underline{v}_j \leq & v_{j,t} \leq \bar{v}_j, \qquad & \forall j,t, \qquad   \label{EqVoltageLimits} \\
& \text{\underline{Ampacity Limits}: } & \quad & l_{j,t} \leq \bar{l}_{j}, \quad & \forall j,t. \qquad \label{EqCurrentLimits}
\end{align}

Constraints \eqref{EqVoltageLimits} and \eqref{EqCurrentLimits} represent voltage and ampacity limits, where $\underline{v}_j$, $\bar{v}_j$, and $\bar{l}_{j}$ are the lower voltage, upper voltage, and line ampacity limits (squared), respectively.
{\color{black} We note that 
the employed branch flow model assumes away the shunt elements 
of the equivalent two-port $\Pi$ line model;
this simplification --- which is mentioned in \cite{FarivarLow2013} --- 
may create some discrepancy in the accuracy of variable $l_{j,t}$ with respect to the ampacity limit $\bar{l}_{j}$,
particularly for underground cables,
as discussed in \cite{NickEtAl_2018}, among others. 
{\color{black}The present work includes the ampacity constraint \eqref{EqCurrentLimits},
but does not delve into the aforementioned potential discrepancy issue, 
since the presented formulation refers to an operational planning optimization problem 
(i.e., scheduling for the next day with an hourly granularity), 
whose limits in both voltages and currents allow sufficient headroom for the real-time operation.}\footnote{{\color{black}
Nevertheless, one of the advantages of the proposed method is the fact that it can dynamically adjust such limits to cater for potential model inaccuracies (see the discussion at the end of Subsection \ref{NonExact}).}
}
}

\underline{Transformer Constraints}, $j \in \mathcal{N}^+_D, t \in \mathcal{T}^+, m = 1,\dots,M$:
\begin{equation} \label{Xf1}
D_{j,t} \geq \alpha_{m} h_{j,t} + \beta_{j, m} l_{j,t} + \gamma_{j,m} \quad \forall j, t, m,
\end{equation}
\begin{equation} \label{Xf2}
h_{j,t} = \delta h_{j,t-1} + \epsilon_j l_{j,t} + \zeta_{j,t},\quad \forall j, t.
\end{equation}
\begin{equation} \label{Cycle}
h_{j,T} = h_{j,0}, \qquad \forall j, 
\end{equation}
Constraints \eqref{Xf1} represent a piecewise linearization 
of the transformer exponential aging acceleration factor
that measures the transformer LoL, 
where $M$ is the number of segments --- see \cite{IEEEtsg20}, 
and the transformer top-oil temperature $h_{j,t}$ is defined by the linear recursive equation \eqref{Xf2}.
The coefficients $\alpha_m$, 
$\beta_{j,m}$,
$\gamma_{j,m}$, 
$\epsilon_{j}$, 
and $\zeta_{j,t}$ are transformer specific,
whereas $\zeta_{j,t}$ also depends on the ambient temperature; 
their detailed formulas and recommended values are found in \cite{IEEEtsg20}. 
Coefficient $\delta$ depends on the time granularity of the problem; 
it is equal to 0.75 for an hourly time period 
and plays an important role in the transformer top oil temperature dynamics. 
Constraint \eqref{Cycle} essentially models the daily 24-hour ahead problem 
as a cycle repeating over an assumed identical next day; 
its dual variable, $\rho_j$, 
captures the impact on the transformer LoL during a future cycle 
resulting from its loading during end of the day hours.

\underline{EV Constraints}, $e \in \mathcal{E}$:
\begin{equation} \label{EVCon2}
\sum_{t \in \mathcal{T}_e} p_{e,t} = \Delta u_{e}, \quad \forall e,
\end{equation}
\begin{equation} \label{EVCon5}
p_{e,t} \leq C_{r}, \quad p_{e,t}^2 + q_{e,t}^2 \leq C_e^2, \quad \forall e, t \in \mathcal{T}_e,
\end{equation}
\begin{equation} \label{EVCon6}
p_{e,t} = q_{e,t} = 0, \quad \forall e, t \in {\mathcal{T^+}} \setminus \mathcal{T}_{e},
\end{equation}
EV constraints \eqref{EVCon2} define the charging requirements, $\Delta u_e$.
Constraints \eqref{EVCon5} impose the charging rate limit (i.e., the EV charger capacity), $C_{r}$, 
and the apparent power limit (i.e., the inverter capacity) $C_e$.
Constraints \eqref{EVCon6} impose zero consumption when the EV is not plugged in, 
with $\mathcal{T}_{e} \subset \mathcal{T^+}$ the set of plugged in time periods.\footnote{
{\color{black} We note that 
since the purpose of this paper is to present an operational-planning formulation for the next day,
we model a stylized formulation of the EV constraints, 
which conveys the EV capabilities and characteristics.
}}

\underline{PV Constraints}, $s \in \mathcal{S}$:
\begin{equation} \label{EqPVcon1}
p_{s,t} \leq \tilde{C}_{s,t}, \quad p_{s,t}^2 + q_{s,t}^2 \leq C_s^2, \quad \forall s, t \in \mathcal{T}_I,
\end{equation}
\begin{equation} \label{EqPVcon2}
p_{s,t} = q_{s,t} = 0, \qquad \forall s, t \in \mathcal{T}^+ \setminus \mathcal{T}_I,
\end{equation}
PV constraints \eqref{EqPVcon1} impose limits on the real and apparent power, 
reflecting the PV nameplate capacity, $C_s$, 
scaled by the intensity of solar irradation $\tilde \rho_t \in [0,1]$ as $\tilde{C}_{s,t} = \tilde \rho_t C_s$. 
Constraints \eqref{EqPVcon2} impose zero generation when  $\tilde \rho_t = 0$, 
where $\mathcal{T}_I \subset \mathcal{T^+}$ the subset of time periods for which $\tilde \rho_t > 0$.

Values in parentheses represent dual variables of the respective constraints.
Variables $l_{j,t}, D_{j,t}, p_{s,t}, p_{e,t} \ge 0$, $\forall j, t, s, e$.

\section{Optimal Grid-DER Coordination} \label{Coordination}

In this section,
we present the proposed framework for achieving the optimal Grid-DER coordination.
In Subsection \ref{Hierarchical}, 
we provide the main idea of our 
decomposition {\color{black} method} to the operational planning problem, 
and we highlight the key differences compared with dual decomposition and ADMM.
Subsection \ref{NetDERopt}, 
describes the network optimization problem and the DER optimization problems 
that constitute the main modules of the proposed {\color{black} method}.
Lastly, in Subsection \ref{NonExact}, 
we provide a ``spacer'' method for dealing with non-exact AC OPF relaxations.

\subsection{Proposed Decomposition Method} \label{Hierarchical}

For illustration purposes, we use an abstract notation, 
which we relate to the enhanced AC OPF problem provided in Section \ref{OPF}.
Let us denote by $\boldsymbol{x}$ the network-related variables (e.g., power flows, voltages, etc.), 
and by $\boldsymbol{y}$ the DER-related variables (e.g., DER real/reactive power generation/consumption).
Also, let $F(\boldsymbol{x})$ represent the network cost included in the objective function \eqref{Obj2},
and $G(\boldsymbol{y})$ the DER cost;
although $G(\boldsymbol{y})$ is zero for the EV/PV models presented in this work, 
we will keep it throughout our analysis since an extension to account for  EV/PV costs is straightforward.
The operational planning problem can then be represented as follows:
\begin{equation} \label{SocialCost}
    \underset{\boldsymbol{x},\boldsymbol{y}}{\text{min}} \quad F(\boldsymbol{x}) + G(\boldsymbol{y}),
\end{equation}
\begin{equation} \label{Balance}
\textit{subject to}:   \quad \boldsymbol{A} \boldsymbol{x} + \boldsymbol{B} \boldsymbol{y} = \boldsymbol{d}, \qquad \qquad
\end{equation}
\begin{equation} \label{NetCon}
    \boldsymbol{x} \in \boldsymbol{\mathcal{X}},
\end{equation}
\begin{equation} \label{DERCon}
    \boldsymbol{y} \in \boldsymbol{\mathcal{Y}},
\end{equation}
where constraint \eqref{Balance} represents the real and reactive power balance constraints \eqref{EqRealBalance}--\eqref{EqReactiveBalance},
containing both network and DER variables, 
$\boldsymbol{\mathcal{X}}$ denotes the set of the remaining networks constraints \eqref{EqVoltageDef}--\eqref{Cycle}, 
and $\boldsymbol{\mathcal{Y}}$ the set of DER constraints \eqref{EVCon2}--\eqref{EqPVcon2}.

{\color{black}
The motive for the coordination between the Grid Operator and the DER owners 
is to maximize social welfare, 
or equivalently minimize the system cost. 
This is the overarching goal of a ``social planner,'' 
in the sense the term is used in welfare economics. 
The coordination problem is written in the form of the centralized AC OPF operational planning problem \eqref{SocialCost}--\eqref{DERCon}, 
which includes both network costs and DER costs/preferences/constraints (as applicable). 
Problem \eqref{SocialCost}--\eqref{DERCon} yields an efficient solution for the entire system (the society), 
i.e., the solution that a ``social planner'' would enforce to maximize social welfare. 
In this solution, DER schedules are optimal for the entire system (society). 

This work does not require a specific relationship between the Grid Operator and the DER owners, 
i.e., it is not required that the power grid utility owns the DERs. 
The basic objective of the Grid Operator is to serve the load, 
and operate the network in a cost-efficient manner, 
while providing DERs with spatiotemporal marginal cost information so that DERs may adapt to the system marginal cost and the joint DER-System actions optimize social surplus. 
This is for example the business model of a municipal distribution utility, 
which has an incentive to drive coordination of DERs and obtain a system-wide-efficient solution. 
The implementation of our method offers the discovery of the marginal cost, 
equivalently the marginal value (e.g., for a rooftop solar) that DERs offer. 
But regardless of the business model, 
which is an open issue related to the ongoing discussions 
on the role of the Grid Operator (or Distribution System Operator), 
the desired outcome of any market design or tariff scheme 
should always be a system wide efficient solution, 
i.e., the solution of \eqref{SocialCost}--\eqref{DERCon}.
}

The key idea presented in this paper is motivated by a pricing interpretation of the operational planning problem {\color{black}and a Grid Operator who aims at finding a system-wide efficient solution}.
Suppose that for some initial forecast of the aggregate{\color{black}d} load, 
{\color{black}the Grid Operator} 
announces to DERs granular prices 
for each location (node) and time period.
DERs, in turn, 
optimize their ``imputed'' income, 
considering that they would be charged/remunerated 
for the real/reactive power consumption/provision 
at the announced spatiotemporal prices.
The optimal DER schedules for the announced prices are then communicated to 
the social planner, 
{\color{black}the Grid Operator}, 
and new prices are calculated and announced to the DERs. 
These new prices are obtained by the dual variables, $\boldsymbol{\lambda}$, 
of the following problem:
\begin{equation} \label{SocialCost2}
    \underset{\boldsymbol{x}}{\text{min}} \quad F(\boldsymbol{x}) ,
\end{equation}
\begin{equation} \label{Balance2}
 \textit{subject to}:   \quad   \boldsymbol{A} \boldsymbol{x} = \boldsymbol{d} - \boldsymbol{B} \boldsymbol{\bar{y}}, \qquad (\boldsymbol{\lambda}) \qquad
\end{equation}
\begin{equation} \label{NetCon2}
    \boldsymbol{x} \in \boldsymbol{\mathcal{X}},
\end{equation}
where $\boldsymbol{\bar{y}}$ represents fixed DER variables (net injections).
The optimization problem \eqref{SocialCost2}--\eqref{NetCon2} 
represents the enhanced AC OPF problem conditional upon the DER dispatch, 
i.e., with known (fixed) net demand.
Hence, 
the dual variables, $\boldsymbol{\lambda}$, 
represent the spatiotemporal marginal costs 
for the specific operating point of the network.
The DERs can then self-schedule in response to the new marginal costs and communicate their new schedules to 
{\color{black}the Grid Operator} 
with the cycle repeated until convergence to a fixed point (up to some tolerance).

To explore tractability and convergence, 
we return to problem \eqref{SocialCost}--\eqref{DERCon}, 
and eliminate network variables $\boldsymbol{x}$, 
by expressing their optimal values as a function of DER variables $\boldsymbol{y}$.
Consider the following problem, for a fixed value of $\boldsymbol{x}$:
\begin{equation} \label{SocialCost3}
    \underset{\boldsymbol{y}}{\text{min}} \quad H(\boldsymbol{d} - \boldsymbol{B y}) + G(\boldsymbol{y}),
\end{equation}
\begin{equation} \label{DERCon2}
\textit{subject to}:  \quad  \boldsymbol{y} \in \boldsymbol{\mathcal{Y}}, \qquad \qquad
\end{equation}
\begin{equation} \label{SocialCost4}
 \qquad  \text{where } \quad   H(\boldsymbol{u}) =  \underset{\boldsymbol{x} \in \boldsymbol{\mathcal{X}, \, \boldsymbol{A x} = \boldsymbol{u}}}{\text{inf}} F(\boldsymbol{x}). 
\end{equation}
If $\boldsymbol{y^*}$ is an optimal solution of \eqref{SocialCost3}--\eqref{SocialCost4}, 
and $\boldsymbol{x^*}$ is an optimal solution of \eqref{SocialCost2}--\eqref{NetCon2}, 
then ($\boldsymbol{x^*}$, $\boldsymbol{y^*}$) is an optimal solution of \eqref{SocialCost}--\eqref{DERCon}.
Following \cite[Prop. 7.2.1]{Ber}, the dual variables, $-\boldsymbol{\lambda}$, obtained by the solution of problem \eqref{SocialCost2}--\eqref{NetCon2} correspond to (sub)gradients of function $H$ at $\boldsymbol{d} - \boldsymbol{B y}$.
Hence, problem \eqref{SocialCost}--\eqref{DERCon} can be solved using an iterative method that updates variables $\boldsymbol{y}$ based on (sub)gradients obtained by solving problems of the form \eqref{SocialCost2}--\eqref{NetCon2}.

The proposed method, 
which falls in the class of proximal algorithms, 
updates the DER variables,  
$\boldsymbol{y}^{(k+1)}$, 
at iteration $k+1$, 
as follows:
\begin{equation} \label{DERupdate}
  \boldsymbol{y}^{(k+1)} =  \text{arg}\underset{{\boldsymbol{y} \in \boldsymbol{\mathcal{Y}}}}{\text{min}} \big\{ G(\boldsymbol{y}) +{\boldsymbol{\lambda}^{(k) T}} \boldsymbol{B y} + \frac{1}{2\sigma^{(k)}} \| \boldsymbol{y} - \boldsymbol{y}^{(k)} \|^2 \big\},
\end{equation}
where $\boldsymbol{\lambda}^{(k)}$ is obtained by solving \eqref{SocialCost2}--\eqref{NetCon2} for $\boldsymbol{\bar y} = \boldsymbol{y}^{(k)}$,
and $\sigma^{(k)}$ is a positive scalar parameter.

Essentially, \eqref{DERupdate} represents a proximal gradient method 
for solving \eqref{SocialCost3}--\eqref{SocialCost4}, i.e.,
\begin{align} \label{DERupdate2}
\begin{split}
\boldsymbol{y}^{(k+1)} & = \text{arg}\underset{\boldsymbol{y} \in \boldsymbol{\mathcal{Y}}}{\text{min}} \big\{ G(\boldsymbol{y})   + \nabla H(\boldsymbol{d} - \boldsymbol{B y})  (\boldsymbol{y} - \boldsymbol{y}^{(k)} )  \\ 
 &   \qquad  \qquad \qquad \qquad + \frac{1}{2\sigma^{(k)}}  \| \boldsymbol{y}  - \boldsymbol{y}^{(k)} \|^2 \big\}. 
  \end{split}
\end{align}
Using $\nabla H(\boldsymbol{d} - \boldsymbol{B y}) = \boldsymbol{B}^T \boldsymbol{\lambda}^{(k)}$,
we have 
\begin{equation*}
 \nabla H(\boldsymbol{d} - \boldsymbol{B y})(\boldsymbol{y} - \boldsymbol{y}^{(k)} )
            = \boldsymbol{\lambda}^{(k)T} \boldsymbol{B} \boldsymbol{y} - \boldsymbol{\lambda}^{(k)T} \boldsymbol{B} \boldsymbol{y}^{(k)},
\end{equation*}
and dropping the constant term $\boldsymbol{\lambda}^{(k)T} \boldsymbol{B} \boldsymbol{y}^{(k)}$
\eqref{DERupdate2} yields \eqref{DERupdate}.
For a sufficiently small positive scalar $\sigma$, 
the method exhibits a cost descent property \cite{Ber}. 
When $G(\boldsymbol{y}) = 0$, it is equivalent to gradient projection (with a stepsize equal to 1),
i.e., 
\begin{equation*}
  \boldsymbol{y}^{(k+1)} =  \Big[ \boldsymbol{y}^{(k)} - \sigma^{(k)} \nabla H(\boldsymbol{d} - \boldsymbol{B y})  \Big]^+,
\end{equation*}
where $[\cdot]^+$ denotes the projection on $\boldsymbol{\mathcal{Y}}$.
Hence, by taking a step $- \sigma^{(k)} \nabla H(\boldsymbol{d} - \boldsymbol{B y})$ 
along the negative gradient, 
and then projecting the result $\boldsymbol{y}^{(k)} - \sigma^{(k)} \nabla H(\boldsymbol{d} - \boldsymbol{B y})$ on $\boldsymbol{\mathcal{Y}}$,
we get (with a stepsize equal to 1) the feasible vector $\boldsymbol{y}^{(k+1)}$.

{\color{black} For clarity, we summarize the proposed 
decomposition method below.

Consider initial schedules $\boldsymbol{y}^{(0)}$.
The proposed method iterates between the following two steps:
\begin{itemize}
    \item \textbf{Step 1}: Solve problem \eqref{SocialCost2}--\eqref{NetCon2} for $\boldsymbol{\bar y} = \boldsymbol{y}^{(k)}$, to obtain $\boldsymbol{\lambda}^{(k)}$.
    \item \textbf{Step 2}: Solve problems \eqref{DERupdate}, to obtain $\boldsymbol{y}^{(k+1)}$.
\end{itemize}
}

Our {\color{black} method} iterates between the network {\color{black}optimization problem} and the DER {\color{black} optimization problems}.
Given a price vector $\boldsymbol{\lambda}^{(k)}$,
DERs respond by setting their quantities (or schedules), $\boldsymbol{y}^{(k+1)}$, according to \eqref{DERupdate}, 
to values that minimize their total cost,
which includes the cost at the given price vector.
The network minimizes the total system cost, according to \eqref{SocialCost2}--\eqref{NetCon2}, 
given the tentatively optimal DER quantities,
i.e., for $\boldsymbol{\bar{y}} = \boldsymbol{y}^{(k+1)}$,
and derives a new price vector, $\boldsymbol{\lambda}^{(k+1)}$,  
reflecting the spatiotemporal marginal costs for the new DER schedules.
Our 
{\color{black} method} involves the iterative solution 
of AC OPF network optimization problems (by network operators)
that are relatively small and simple 
since they are conditional upon DER schedules, 
followed by parallelizable individual DER self-scheduling 
responding to tentative spatiotemporal prices/DLMCs provided by the network problem solution.
At the same time, 
the information exchanged between network operators and DERs 
is limited to DLMCs in one direction and real/reactive power services in the other. 
This renders network information 
and DER preference and capability information 
internal to the network operator solving the {\color{black}network optimization} problem 
and the DERs solving the self-scheduling {\color{black} optimization} problems in parallel 
and not necessarily in a synchronous lock step fashion.
We note that the magnitude of regularization terms included in the DER {\color{black} optimization} problems 
represents the willingness of each DER to accept a sub-optimal solution at each iteration 
limiting the distance of their new schedule from the previous iteration's schedule. 
This is precisely how synchronization oscillations are attenuated, 
and singularities (e.g., from instances of DLMC equal to zero) are handled. 
Moreover, proper scaling of synchronization terms may speed up convergence.

{\color{black}We wish to emphasize that there are key differences of our {\color{black} method}, 
which falls in the class of proximal algorithms,
compared to dual decomposition and ADMM.}
We thus highlight them here.
Both dual decomposition and ADMM (see also their summary provided in \cite{MolzahnEtAl_2017}) 
employ a Lagrangian relaxation (augmented in the case of ADMM),
and include steps that minimize the Lagrangian
followed by updates of the dual variables, $\boldsymbol{\lambda}$,
driven by relaxed constraint imbalances.
These dual variables are then provided to the network and the DERs to update the primal variables.
Consider, for instance, an iteration of dual decomposition:
\begin{equation} \label{SocialCostDual}
\boldsymbol{x}^{(k+1)} = \text{arg}\underset{\boldsymbol{x} \in \boldsymbol{\mathcal{X}}}{\text{min}} \, \{ F(\boldsymbol{x}) + \boldsymbol{\lambda}^{(k)T} \boldsymbol{A} \boldsymbol{x} \},
\end{equation}
\begin{equation} \label{DERCostDual}
 \boldsymbol{y}^{(k+1)} = \text{arg}\underset{\boldsymbol{y} \in \boldsymbol{\mathcal{Y}}}{\text{min}} \, \{ G(\boldsymbol{y}) + \boldsymbol{\lambda}^{(k)T} \boldsymbol{B} \boldsymbol{y} \},
\end{equation}
\begin{equation} \label{lambda-update}
  \boldsymbol{\lambda}^{(k+1)} = \boldsymbol{\lambda}^{(k)} + \sigma^{(k+1)} ( \boldsymbol{A} \boldsymbol{x}^{(k+1)} + \boldsymbol{B} \boldsymbol{y}^{(k+1)} - \boldsymbol{d}).   
\end{equation}
Essentially, 
problem \eqref{SocialCostDual} selects the flows 
by choosing some ``generic'' DERs at each node at a cost equal to $\lambda^{(k)}$.\footnote{This is in fact reminiscent of the methodology employed in \cite{LMV},
where we assumed ``generic'' DERs in the context of non wires alternatives.}
Evidently, problem \eqref{DERCostDual} updates the DER variables, 
similarly to \eqref{DERupdate} without the augmentation terms.
Notably, an appropriate ADMM implementation would include these terms and update DER variables similarly to \eqref{DERupdate}.
{\color{black} Despite the similarity in the $y-$update between \eqref{DERCostDual} and \eqref{DERupdate} 
--- in other words between the $y-$update of ADMM (or dual decomposition) and Step 2 of the proposed method, 
there is a fundamental difference between the $x-$update and $\lambda-$update of ADMM (or dual decomposition) 
--- see \eqref{SocialCostDual} and \eqref{lambda-update} --- 
and Step 1 of the proposed method.  
Both ADMM and dual decomposition perform the $\lambda-$update 
using the imbalance quantity $\boldsymbol{A} \boldsymbol{x}^{(k+1)} + \boldsymbol{B} \boldsymbol{y}^{(k+1)} - \boldsymbol{d}$,
whereas the proposed method calculates gradient information, i.e., $\boldsymbol{\lambda}$ from the solution of \eqref{SocialCost2}--\eqref{NetCon2}.
Hence, the proposed method, unlike dual decomposition and ADMM, enjoys convergence properties
of proximal algorithms even in non-convex settings.
Most importantly, even under convexity assumptions, }
a common drawback of both dual decomposition and ADMM 
is that until convergence is reached,
the solution obtained is not feasible,
i.e., the power balance equations are not satisfied for the obtained DER variables.
On the contrary, 
our {\color{black} method} yields a solution that satisfies the power balance at each iteration.
From a practical perspective, this is a crucial point, 
{\color{black}and major advantage of the proposed method.}

\subsection{Network and DER Optimization Problems} \label{NetDERopt}

This subsection formulates the details of the network and DER optimization problems discussed above in   
the higher level abstract notation. 
The network-related variables $\boldsymbol{x}$ include 
power flows $P_{j,t}$, $Q_{j,t}$, 
voltages $v_{j,t}$, 
currents $l_{j,t}$, 
and transformer variables $D_{j,t}$ and $h_{j,t}$.
DER variables include 
EV/PV real/reactive power schedules $p_{e,t}$, $q_{e,t}$, $p_{s,t}$, and $q_{s,t}$.

\subsubsection{Network {\color{black}Optimization} Problem}
For given EV/PV schedules, 
$p_{e,t}^{(k)}$, $q_{e,t}^{(k)}$, $p_{s,t}^{(k)}$, and $q_{s,t}^{(k)}$,
the power balance constraints are:
\begin{align} \label{EqRealBalance2}
\begin{split}
P_{j,t} & - r_{j} l_{j,t}  = \sum_{j' \in \mathcal{C}_j} P_{j',t} + p^d_{j,t} \\ 
& + \sum_{e \in \mathcal{E}_{j,t}}  p^{(k)}_{e,t} - \sum_{s \in \mathcal{S}_j} p^{(k)}_{s,t} \rightarrow (\lambda^{P \,(k)}_{j,t}), \quad \forall j,t,
\end{split}
\end{align}
\begin{align} \label{EqReactiveBalance2}
\begin{split}
Q_{j,t} & - x_{j} l_{j,t}  = \sum_{j' \in \mathcal{C}_j} Q_{j',t} + q^d_{j,t} \\ 
& + \sum_{e \in \mathcal{E}_{j,t}}  q^{(k)}_{e,t} - \sum_{s \in \mathcal{S}_j} q^{(k)}_{s,t} \rightarrow (\lambda^{Q \, (k)}_{j,t}), \quad \forall j,t.
\end{split}
\end{align}
Hence, the fixed injections network {\color{black} optimization} problem, corresponding to \eqref{SocialCost2}--\eqref{NetCon2} is:
\begin{equation*}
 \textbf{Net-opt: } \text{minimize } \eqref{Obj2}
\textit{ subject to}: \,\, \eqref{EqRealBalance2}-\eqref{EqReactiveBalance2}, \eqref{EqVoltageDef}-\eqref{Cycle},   
\end{equation*}
and its solution yields dual variables $\lambda^{P \,(k)}_{j,t}$ and $\lambda^{Q \,(k)}_{j,t}$.
However, since fixed injections may occasionally cause  voltage and/or ampacity limits to be violated, we replace hard constraints by soft  constraints introducing  
\begin{equation} \label{VoltageSoft}
 \Delta v_j \geq v_j -\bar{v}_j, \, \Delta v_j \geq \underline{v}_j - v_j, \, \Delta v_j \geq 0, \, \forall j \in \mathcal{N}^+,
\end{equation}
\begin{equation} \label{CurrentSoft}
 \Delta l_j \geq l_j -\bar{l}_j, \, \Delta l_j \geq 0, \, \forall j \in \mathcal{N}^+,
\end{equation}
and adding penalty terms $M^v \sum_{j} (\Delta v_j)^2$, $M^l \sum_{j} (\Delta l_j)^2$, to the objective function,
with penalty coefficients $M^v$, $M^l$.
{\color{black}The reason for selecting a quadratic form is to obtain a linear first derivative, (which affects the marginal cost), instead of a constant first derivative, which would have been the case, had we selected linear penalties. Since our objective is to derive the marginal costs from the {\color{black} network optimization} problem, a quadratic penalty captures more accurately the value of the hard constraint's dual by selecting a slack at which the derivative of the penalty term equals to or is indeed very close to the marginal cost of enforcing the constraint. In fact, this softer approach facilitates the 
decomposition's convergence.}{\footnote{\color{black}
We also note that the ampacity limits,
similarly to the transformer limits, 
may themselves be in reality, not hard constraints 
--- e.g., the cable limits also depend on the thermal degradation; 
in fact, in our parallel work \cite{CIRED}, 
we are dealing with this issue and illustrate that a similar to the transformer model can be employed to account for cable aging.}}

\subsubsection{DER {\color{black} Optimization} Problems}
Applying \eqref{DERupdate}, we obtain the following PV/EV self-scheduling optimization problems:
\paragraph{PV {\color{black} Optimization} Problem}
Each self-scheduling PV, $s$, 
adapts its real/reactive power profile, 
$p_{s,t}$, $q_{s,t}$, 
to the tentative location specific real and reactive power prices, 
$\lambda_t^{P\,(k)}$, $\lambda_t^{Q \, (k)}$, 
by selecting the power factor of its smart inverter to  maximize revenues from the provision of real and reactive power over the optimization horizon. Namely it solves:
\begin{equation} \label{PVmax1}
\begin{split}
\textbf{PV-opt}: \qquad \,\,\, \underset{p_{s,t},q_{s,t}}{\text{max}}  
\sum_t { \lambda_t^{P\,(k)} p_{s,t} } 
+ \sum_t { \lambda_t^{Q\,(k)} q_{s,t}}    \\
- \frac{1}{2\sigma^{(k)}} \sum_t { \big( p_{s,t} - p_{s,t}^{(k)} \big) ^2}
- \frac{1}{2\sigma^{(k)}}  \sum_t { \big( q_{s,t} - q_{s,t}^{(k)} \big) ^2},
\end{split}
\end{equation}
\emph{subject to}: PV constraints \eqref{EqPVcon1}-\eqref{EqPVcon2}.
Note that PV-opt is separable across time periods.
If $\lambda_t^{P\,(k)}$ is negative, 
then PV $s$ sets $p_{s,t} = 0$;
If $\lambda_t^{Q\,(k)}$ is positive (negative), 
then PV $s$ adjusts its power factor to provide (consume) reactive power, 
i.e., $q_{s,t} > 0$ ($q_{s,t} < 0$), 
to render  $\lambda_t^{Q\,(k)} q_{s,t}$ positive.

\paragraph{EV {\color{black} Optimization} Problem}
Each EV, $e$, 
solves the EV-opt problem shown below to  adapt its real/reactive power charging profile, 
$p_{e,t}$, $q_{e,t}$, 
to the current iteration's real and reactive power prices, 
$\lambda_t^{P\,(k)}$, $\lambda_t^{Q \, (k)}$, 
so as to minimize over the daily cycle its charging cost net of revenues from reactive power services. The EV self-scheduling problem is:
\begin{equation} \label{EVmin2}
\begin{split}
\textbf{EV-opt}: \qquad \,\,\, \underset{p_{e,t},q_{e,t}}{\text{min}}  
\sum_t { \lambda_t^{P\,(k)} p_{e,t} } 
+ \sum_t { \lambda_t^{Q\,(k)} q_{e,t}}  \\
 +  \frac{1}{2\sigma^{(k)}} \sum_t { \big( p_{e,t} - p_{e,t}^{(k)} \big) ^2}
+ \frac{1}{2\sigma^{(k)}} \sum_t { \big( q_{e,t} - q_{e,t}^{(k)} \big) ^2},
\end{split}
\end{equation}
\emph{subject to:} EV inter-temporally coupled constraints \eqref{EVCon2}--\eqref{EVCon6}.
When $\lambda_t^{Q\,(k)}$ is positive (negative), 
then the EV has the opportunity to generate income by providing (consuming) reactive power, 
and if it is plugged in during multiple time periods it chooses to charge during low energy cost hours.

\subsection{Dealing with Non-Exact Relaxations} \label{NonExact}

This subsection considers instances when the SOCP relaxations are not exact, 
in which primal solutions that violate the laws of physics
provide fallacious dual variable information (price signals)
to the DERs.
Consider, e.g., iteration $k$, with
\begin{equation*}
\text{Gap}^{(k)} = \sum_{j,t} ( v^{(k)}_{u(j),t} l^{(k)}_{j,t} - {P^{(k)}_{j,t}}^2 - {Q^{(k)}_{j,t}}^2 ) > \tau,    
\end{equation*}
where $\tau$ is a small tolerance.
Following \cite{WeiEtAl_2017},
we add an additional constraint in the opposite direction to enforce equality 
\begin{equation*}
 v_{u(j),t} l_{j,t} \leq {P_{j,t}}^2 + {Q_{j,t}}^2, \qquad \forall j, t,  
\end{equation*}
which can be rewritten as 
\begin{equation*}
 (v_{u(j),t} + l_{j,t})^2  \leq  (v_{u(j),t} - l_{j,t})^2 + (2P_{j,t})^2 + (2Q_{j,t})^2.
\end{equation*}
The idea in \cite{WeiEtAl_2017} is 
to replace the \emph{rhs} by its linearization around $\boldsymbol{x}^{(k)}$, 
and penalize the gap, $w_{jt,}$, 
with a carefully tuned penalty $\bar \rho^i$, 
where $i$ is an inner loop iteration updating the penalty and solving repeatedly the following problem
till the gap is within tolerance:
\begin{equation} \label{InexactObj}
 \text{min} 
 \sum_{t} \big( c^P_t P_{1,t}
+ c^Q_t Q_{1,t} 
+ \sum_{j \in \mathcal{N}_L} c^D_j D_{j,t} 
+ \bar \rho^{i} \sum_{j} w_{j,t} 
\big),    
\end{equation}
\begin{equation*}
\textit{subject to}: \qquad \qquad \eqref{EqRealBalance2}-\eqref{EqReactiveBalance2}, \eqref{EqVoltageDef}-\eqref{Cycle}, \qquad     \qquad \qquad
\end{equation*}
\begin{align} \label{LinConcave}
\begin{split}
(v_{u(j),t} & + l_{j,t})^2 - \big[  (v^{(k)}_{u(j),t} - l^{(k)}_{j,t})^2 + 4{P^{(k)}_{j,t}}^2 + 4{Q^{(k)}_{j,t}}^2 \\
            & + 2(v^{(k)}_{u(j),t} - l^{(k)}_{j,t})v_{u(j),t} - 2(v^{(k)}_{u(j),t} - l^{(k)}_{j,t})l_{j,t} \\
            & + 8P^{(k)}_{j,t}P_{j,t} + 8Q^{(k)}_{j,t}Q_{j,t} \big] \leq w_{j,t}, \, \forall j, t,      
\end{split}
\end{align}
\begin{equation} \label{LinGap}
    w_{j,t} \geq 0, \qquad \forall j,t.
\end{equation}

Note that, although the above inner loop is used in conjunction with the fixed injections Net-Opt problem, 
it is also applicable to the centralized enhanced AC OPF problem presented in Section \ref{OPF}.
An important point requiring attention 
in order to obtain correct DLMCs that drive DER schedule updates, 
is that, once the operating point of the network that observes the laws of physics is obtained, 
Net-opt of Subsection \ref{NetDERopt} must be solved again as a convex problem after replacing \eqref{EqCurrentDef} by the linearized equality constraint:
\begin{equation} \label{LinCur}
    l_{j,t} = \frac{2 {P^{(k)}_{j,t}}}{v^{(k)}_{u(j),t}} P_{j,t}
    + \frac{2 {Q^{(k)}_{j,t}}}{v^{(k)}_{u(j),t}} Q_{j,t}
    - \frac{ {P^{(k)}_{j,t}}^2 + {Q^{(k)}_{j,t}}^2}{{v^{(k)}_{u(j),t}}^2} v_{u(j),t}, \forall j,t.
\end{equation}

{\color{black}
To deal with instances of non-exact relaxation, 
we expand Step 1 to diagnose non-exact relaxations and to establish the physics of the Load Flow problem.

\begin{itemize}
    \item \textbf{Step 1}: Solve problem \eqref{SocialCost2}--\eqref{NetCon2} for $\boldsymbol{\bar y} = \boldsymbol{y}^{(k)}$.
    \begin{itemize}
    \item \textbf{Step 1.a}: If constraint \eqref{EqCurrentDef} is binding at equality (exact relaxation), then obtain $\boldsymbol{\lambda}^{(k)}$ and proceed to Step 2.
    \item \textbf{Step 1.b:} Set $i = 0$.
    \item \textbf{Step 1.c}: Solve \eqref{InexactObj}, subject to: \eqref{EqRealBalance2}--\eqref{EqReactiveBalance2}, \eqref{EqVoltageDef}--\eqref{Cycle}, \eqref{LinConcave}--\eqref{LinGap}.
    \item \textbf{Step 1.d}: If $\text{Gap}^{(k)} > \tau$, then update $i \leftarrow i+1$ and return to Step 1.c.
    \item \textbf{Step 1.e}: Repeat Step 1.c replacing \eqref{EqCurrentDef} by \eqref{LinCur}, to obtain $\boldsymbol{\lambda}^{(k)}$.
    \end{itemize}
    \item \textbf{Step 2}: Solve problems \eqref{DERupdate}, to obtain $\boldsymbol{y}^{(k+1)}$.
\end{itemize}

The proposed ``remedy'' recovers the solution of the non-convex AC OPF problem,
and retrieves the ``correct'' gradient information.
Notably, as suggested by \cite{WeiEtAl_2017}, the AC OPF solution can be fine-tuned (if required) 
by warm-starting a non-convex solver from the operating point identified by the proposed remedy.

Last but not least, we should note that one of the advantages of the proposed iterative algorithm 
is that we can dynamically adjust limits/parameters at each iteration, 
to deal with model inaccuracies (e.g., the ampacity limits as discussed in Section \ref{OPF}).
A similar process (but only with a few iterations) is actually followed by ISOs to ensure e.g., 
AC-feasibility of the solution or the satisfaction of certain contingencies. 
In practice, the solution is checked, 
and if constraints that are not captured by the formulation are found to be violated, 
these constraints are added to the formulation 
and/or parameters/limits are adjusted to ensure feasibility. 
In day-ahead markets, this process is repeated 
until either a feasible solution is found or the ISO runs out of time.
In the context of our AC OPF problem, the basic advantage of the proposed framework is that 
we can check the solution (using e.g., more accurate Load Flows, and adjusting ampacity limits) 
whenever we observe violations.
This advantage is basically attributed to the fact that our proposed method can ensure feasibility at each iteration.
In order to guarantee convergence, we only need to keep these parameters fixed in the end.
}

\begin{figure*}[tb]
\centering
\includegraphics[width=6.5in]{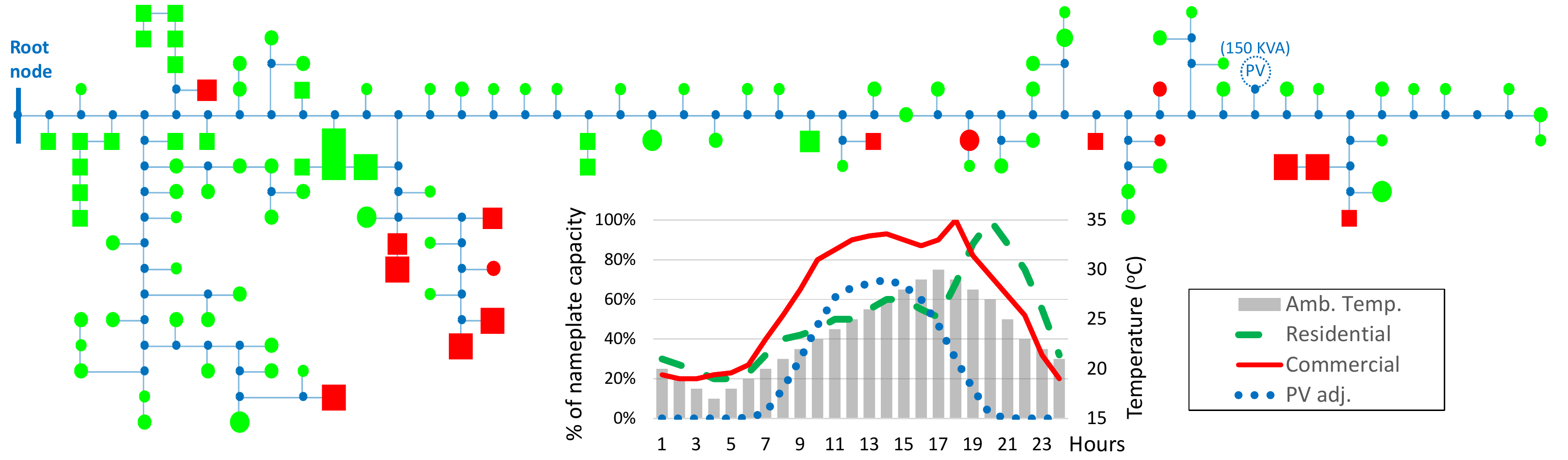}
\caption{Feeder topology diagram (307 nodes). 110 transformers: round (pole); square (pad); green (residential); red (commercial); sizes reflect nameplate capacity. Residential and commercial load profiles, PV adjustment factor ($\tilde \rho_t$), and ambient temperature. PV installation 150 {\color{black}k}VA (p.f. = 1).} 
\label{figFeeder}
\end{figure*}

\section{Test Case} \label{Numerical}
This section 
illustrates the implementation of the proposed 
decomposition method
on an actual 13.8-{\color{black}k}V feeder of Holyoke Gas and Electric (HGE), 
a municipal distribution utility in MA, US.
Subsection \ref{Input} presents briefly the input data.
Subsection \ref{Remarks} provides some computational and practical implementation remarks.
Subsection \ref{Results} presents and discusses the main numerical results.
{\color{black}Subsection \ref{sens} provides additional results and sensitivity analysis.}

\subsection{Input Data} \label{Input}
\begin{table}[t] 
	\caption{Aggregate Line and Transformer Data}  \label{tab1} 
\centering
\includegraphics[width=3.35in]{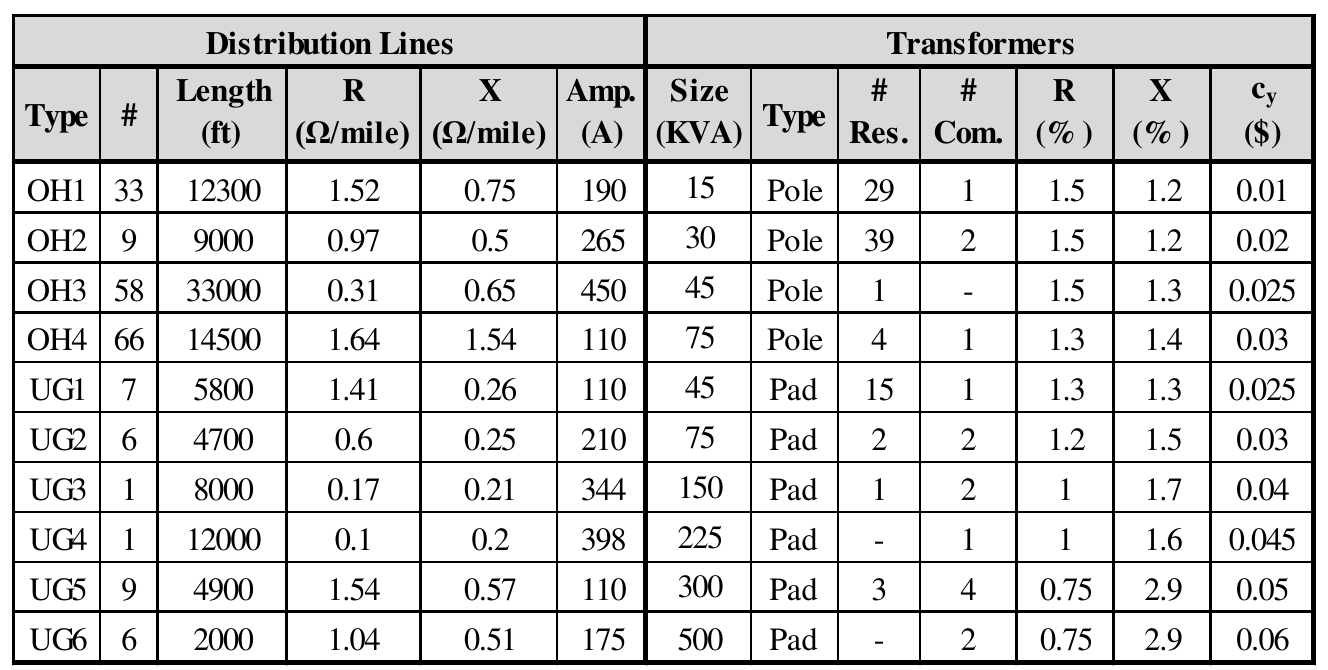}
\end{table}
Feeder topology is presented in Fig. \ref{figFeeder}
that also shows residential and commercial load profiles (obtained from \cite{CIGRE})
as a percentage of transformer nameplate capacity,
PV adjustment factor, $\tilde \rho_t$,
and ambient temperature. 
A power factor of 0.95 (0.85) is assumed for residential (commercial) nodes.
Aggregate line and transformer data are listed in Table \ref{tab1}.
Transformer-specific parameters are obtained by applying the formulas presented in \cite{IEEEtsg20}.
Lower (upper) voltage limits are set equal to 0.95 (1.05) p.u.
LMPs range from 25.59 to 53.48 (\$/MWh).
Reactive power opportunity cost is assumed equal to 10\% the value of the LMP.

We considered an EV scenario, 
elaborating data from \cite{NHTS}. 
In total, we allocated 662 EVs,
with charging requirements that ranged from 5.97 to 47.54 {\color{black}k}Wh, 
connected from 5 to 21 hours during the day, 
with a penetration that ranged from: 
1 to 3 EVs in the 15-{\color{black}k}VA transformers, 
3 to 6 EVs in the 30-{\color{black}k}VA, 
4 to 8 EVs in the 45-{\color{black}k}VA, 
8 to 12 EVs in the 75-{\color{black}k}VA, 
and 12 to 24 EVs in the remaining 150-{\color{black}k}VA (and higher) transformers.
We also considered a PV scenario,
assuming a penetration of two 10-{\color{black}k}VA rooftop solar at each transformer.

\subsection{Computational and Practical Remarks} \label{Remarks}

\subsubsection{Transformer Degradation Piecewise Linearization}

Linearization breakpoints introduce a difficulty in the correct estimation of associated duals,
and may create some inter-iteration oscillations.
To mitigate oscillatory behavior, 
we started by solving the centralized network problem 
with breakpoints at 100, 105, 110, 115, 120, 130, 140, and 150$^o$C, 
proceeded to gradually increase the density around the optimal HST from 5$^o$C to 0.5$^o$C,
and finally kept the high-density linearization.

\subsubsection{Future Impact of Transformer Degradation}

An interesting point taken up in \cite{IEEEtsg20} 
is the impact of scheduling decisions (due to impact of transformer loading) on costs incurred beyond the optimization horizon. 
To ensure consistent comparisons, 
the centralized problem was modeled as a repeating 24-hour cycle 
yielding a dual to constraint \eqref{Cycle}.
In the numerical results reported here,  
we assumed an initial condition of $h_{j,0}$ 
and appended constraint \eqref{Cycle} to the objective function, 
using the dual value obtained as described above.

\subsubsection{DLMC Initialization}
In the reported numerical results, 
we estimated initial DLMCs 
by assuming that EVs and PVs respond to the substation LMPs 
and do not provide reactive power services.
As such, 
initial DLMCs may be viewed as corresponding to 
open loop Time of Use prices \cite{IEEEtsg20}.

\subsubsection{Regularization Terms}
In the first iteration, 
regularization terms were dropped 
allowing EVs/PVs to freely select their real/reactive power schedules 
based on initial DLMCs. 
In subsequent iterations, 
regularization terms were used with $\sigma = 10^{-4}$ 
that proved to be small enough 
to avoid oscillations during early iterations.
To guarantee descent while close to the system optimal solution, 
$\sigma$ was later reduced by a factor of $\frac{2}{3}$.

\subsubsection{Computational Times}
The SOCP optimization problems were solved on a Dell Intel Core i7-5500U @2.4 GHz with 8 GB RAM, 
using CPLEX 12.7. 
Solution times for the SOCP AC OPF problems were in the order of 5 to 20 sec.

\subsection{Numerical Results} \label{Results}

\begin{figure}[t]
\centering
\includegraphics[width=3.45in]{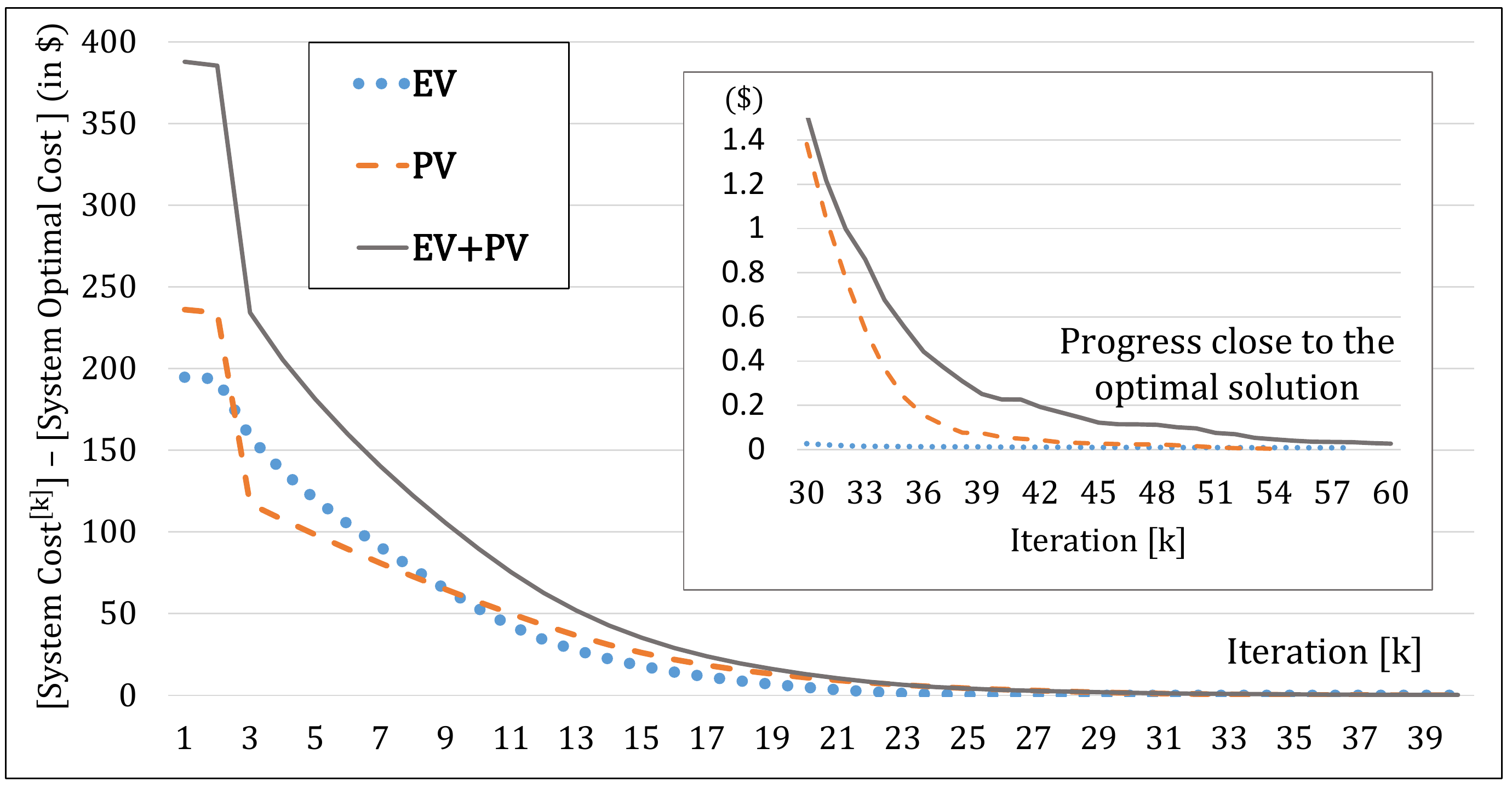}
\caption{Convergence progress of the 
decomposition {\color{black} method}.} 
\label{FigR1}
\end{figure}

We first consider three scenarios: 
the \textbf{EV} scenario (662 EVs), 
the \textbf{PV} scenario (220 PVs), 
and a combined \textbf{EV+PV} scenario (662 EVs and 220 PVs).
Fig. \ref{FigR1} 
illustrates the convergence of the proposed {\color{black} method},
relative to a centralized optimal solution benchmark; 
the latter in absolute numbers (\$) is 3766.46, 2758.51, and 3054.49 
for the \textbf{EV}, \textbf{PV}, and the combined \textbf{EV+PV} scenario, 
respectively.
In Fig. \ref{FigR1},
we show the difference of the system cost at each iteration from the optimal solution,
and observe that convergence (in the order of 0.01\$) is practically achieved 
within a few tens of iterations. 

\begin{figure}[t]
\centering
\includegraphics[width=3.45in]{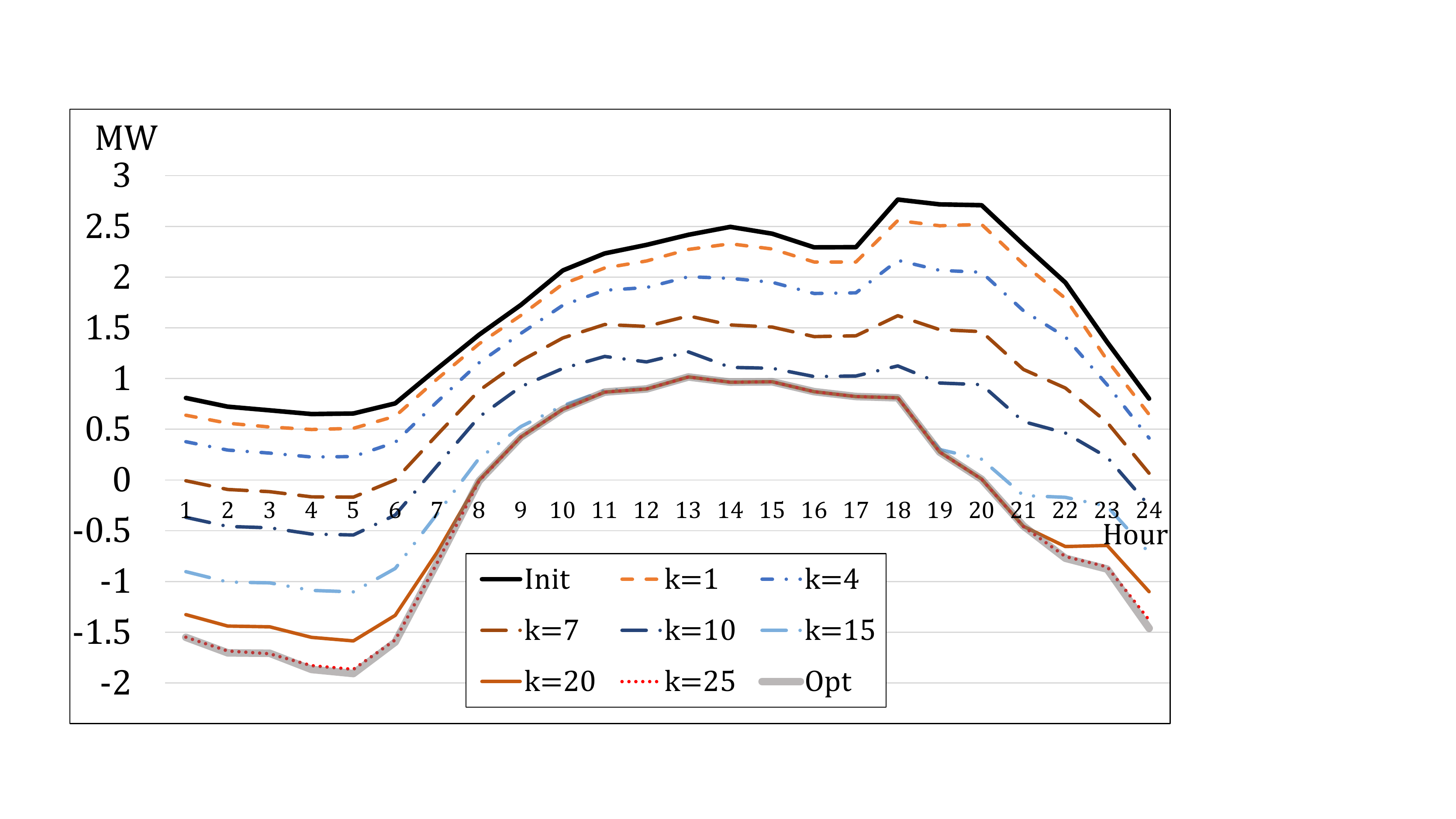}
\caption{Reactive power flow at the root node for iteration $k$; \textbf{EV} scenario. Init: EVs respond to LMP, operating with p.f. = 1. Opt: System optimal solution.} 
\label{FigR2}
\end{figure}
\begin{figure}[t]
\centering
\includegraphics[width=3.45in]{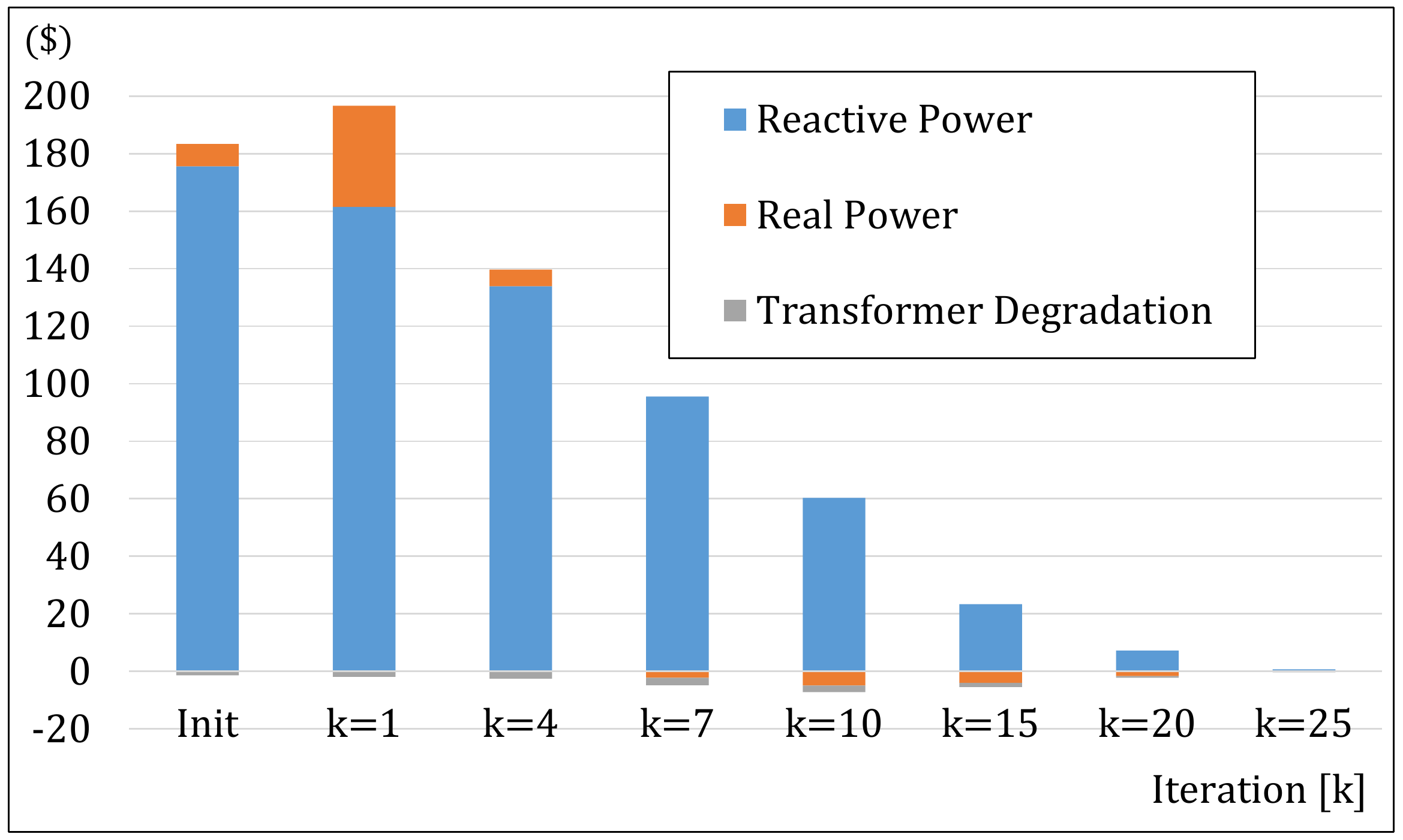}
\caption{System cost deviation relative to optimal solution for iteration $k$, per component: 
reactive power, real power, transformer degradation; \textbf{EV} scenario.} 
\label{FigR3}
\end{figure}

Taking a closer look at the \textbf{EV} scenario,
we present in Fig. \ref{FigR2} the reactive power flow at the substation.
Initially, 
EVs respond to LMP, 
with no reactive power provision (see Init),
and gradually they adjust their reactive power profile,
to reach the system optimal profile (see Opt).
In fact, in the \textbf{EV} scenario,
the provision of reactive power is the major contributor to system cost reduction 
as shown in Fig. \ref{FigR3}.
This graph shows the difference of the three components of the system cost 
(real power, reactive power, transformer degradation) compared to the optimal solution.
It explains the trajectory of Fig. \ref{FigR1}, 
and the fact that the major benefit in the system cost comes from the reactive power provision, 
as shown in Fig. \ref{FigR2}.

\begin{figure}[tb]
\centering
\includegraphics[width=3.45in]{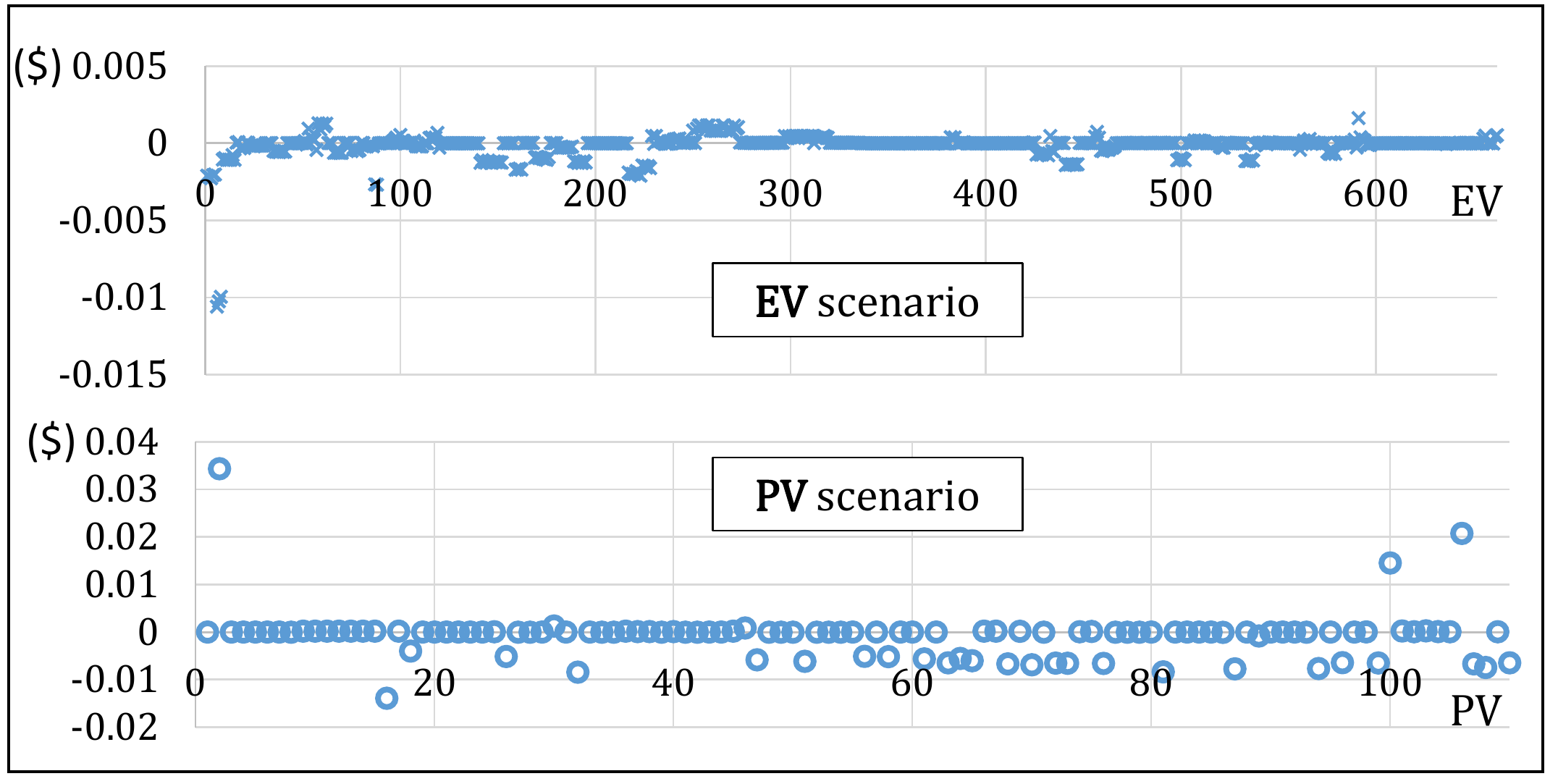}
\caption{Imputed cost (revenue) deviation relative to optimal solution; \textbf{EV} (\textbf{PV}) scenario.} 
\label{FigR4}
\end{figure}

In Fig. \ref{FigR4}, we present the deviation from the optimal of 
the imputed costs for the \textbf{EV} scenario, 
and the imputed revenues for the \textbf{PV} scenario,
for the last iteration presented in Fig. \ref{FigR1}.
Reported deviations are all relative to the costs/revenues calculated by the benchmark centralized solution.
To get a sense of the costs, 
on average the EV imputed cost is 0.48\$, 
whereas the PV imputed revenue is 2.32\$.
The observed differences are very small.
For instance, in the \textbf{PV} scenario, 
the highest difference represents less than 1.5\% of the PV imputed revenue.
Our numerical experiments showed that 
these small differences are mainly due to the transformer degradation linearization, 
resulting, 
if not dense enough, 
in oscillatory behavior impeding convergence of the duals.
Nevertheless, 
the estimated DLMCs were very close to the values in the centralized solution. 
For instance, 
consider the DLMCs at the transformer locations
for the \textbf{EV} scenario:
90\% of P-DLMCs and the Q-DLMCs are with 0.01\$/MWh, 
whereas most of the remaining 10\% are within 0.1\$/MWh.
Only one P-DLMC (Q-DLMC) instance reaches close to 1.5\$/MWh (\$/MVARh).

\begin{figure}[tb]
\centering
\includegraphics[width=3.45in]{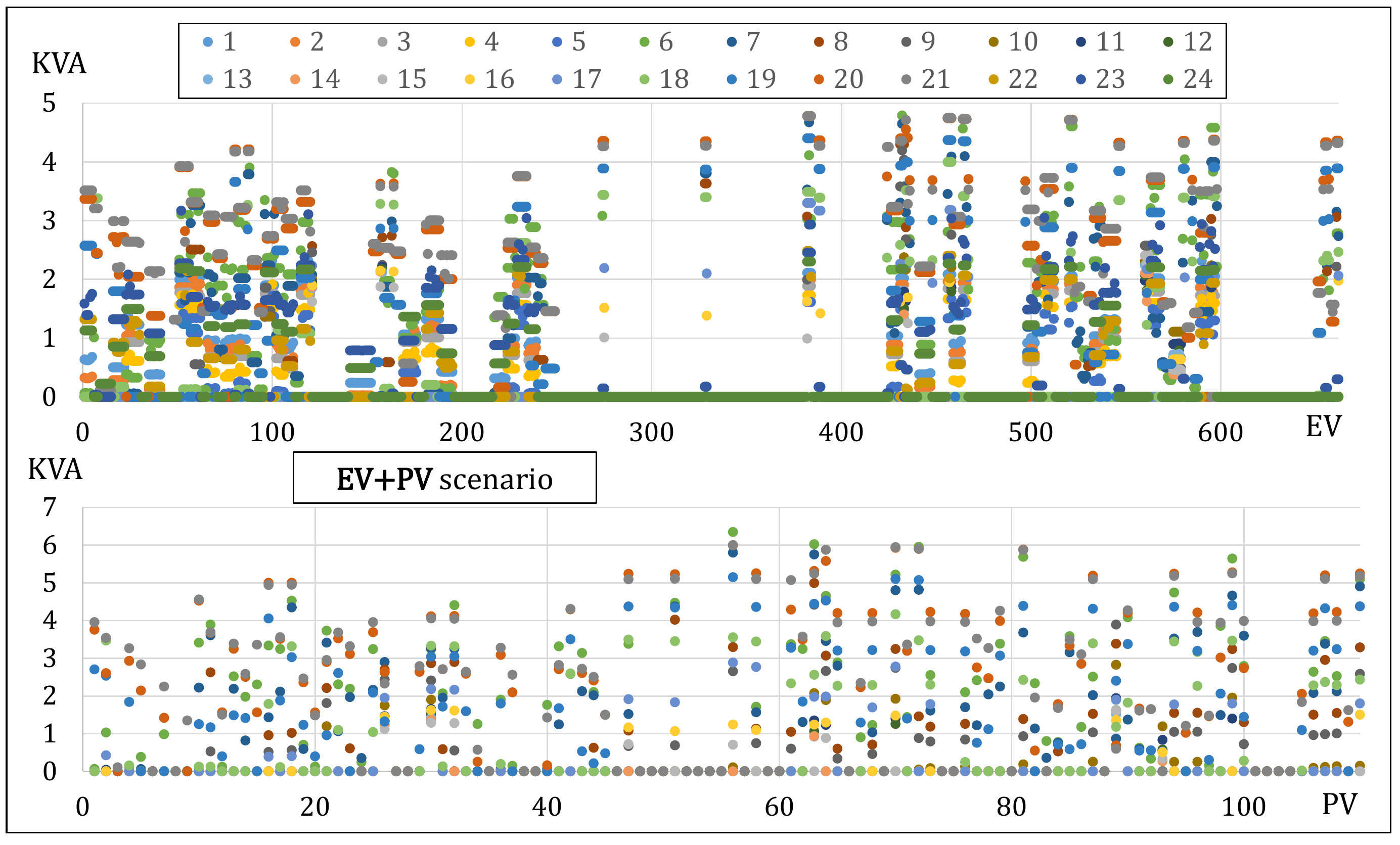}
\caption{Spare inverter capacity; \textbf{EV+PV} scenario.} 
\label{FigR5}
\end{figure}

Another interesting observation is the utilization of the inverter, 
which relates to
the amount of reactive power provided by each resource.
We showed in Fig. \ref{FigR2} that 
the aggregate reactive power provision soon reaches the optimum 
for the \textbf{EV} scenario. 
This is a general remark for all scenarios.
EVs and PVs in fact provide the ``correct'' amount of reactive power 
even with zero or nearly zero DLMCs 
(where by ``correct'' we mean system-optimal as obtained by the centralized solution).
Fig. \ref{FigR5} shows the under-utilization of the inverter, 
i.e., 
the spare inverter capacity of EVs and PVs under the \textbf{EV+PV} scenario.
Unsurprisingly, 
the DLMCs in the centralized solution, 
at the locations and time periods where an under-utilization is observed,
are practically zero.


\subsection{Additional Results and Sensitivity Analysis} \label{sens}

In this subsection, we further elaborate on the aforementioned scenarios,
and provide additional numerical results along with sensitivity analysis as follows: 
1) Increased EV penetration, with emphasis put on the increased transformer degradation;
2) Negative LMPs, which render the SOCP relaxation non-exact, with emphasis on the penalty selection parameters;
3) Sensitivity Analysis of the Reactive Power Opportunity Cost, and
4) Comparison with Time of Use Tariffs.

\subsubsection{{\color{black}Increased EV penetration}}
\begin{figure}[tb]
\centering
\includegraphics[width=3.45in]{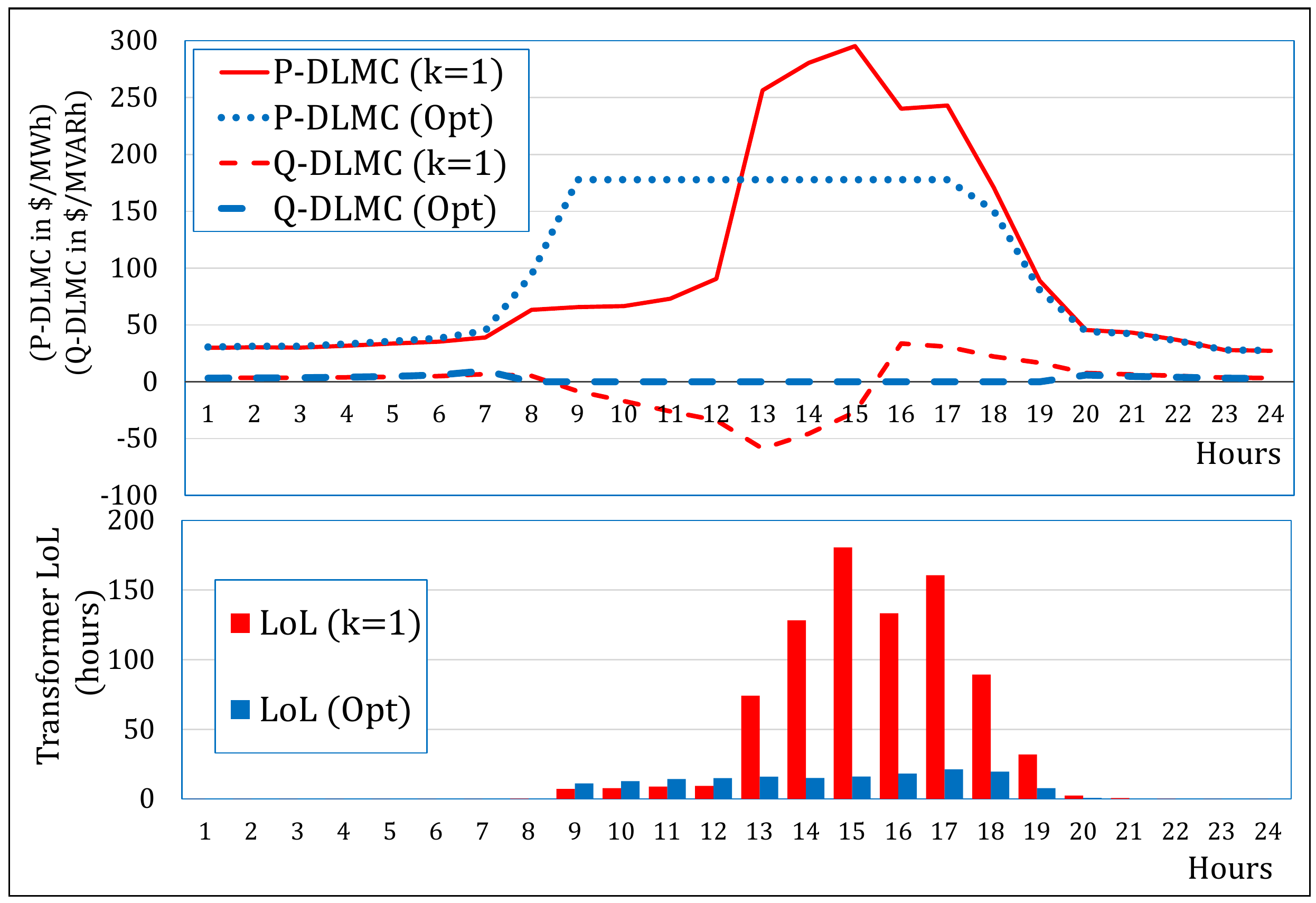}
\caption{P-DLMCs and Q-DLMCs (top) and hourly Loss of Life (bottom) at transformer location with the highest observed P-DLMC, for iteration $k=1$  in comparison to the optimal solution; \textbf{EVx2} scenario.}
\label{FigR7}
\end{figure}

We first tested a high EV penetration scenario, \textbf{EVx2}, 
duplicating the 662 EVs.
\textbf{EVx2} required soft limits \eqref{VoltageSoft}--\eqref{CurrentSoft},
with $M^v = 5000$, and $M^l= 1000$; 
some small voltage violations were observed in 3 nodes.
An interesting remark is the mutual adaptation of DLMCs, 
and transformer hourly LoL converging to smooth profiles,
as shown in Fig. \ref{FigR7},
for the node with the highest observed P-DLMC.
Both P-DLMC and Q-DLMC curves (top) are flattened, 
by converging EV profiles resulting in  
smoother LoL profiles (bottom). 
Notably, 
aggregate LoL is reduced from 835 hours (at $k=1$) to 168 hours (at the optimal solution).
Q-DLMCs also become nearly zero for several hours at the optimal solution,
flattening the positive/negative spikes exhibited at $k=1$.

\subsubsection{{\color{black}} Negative LMPs}
\begin{figure}[tb]
\centering
\includegraphics[width=3.45in]{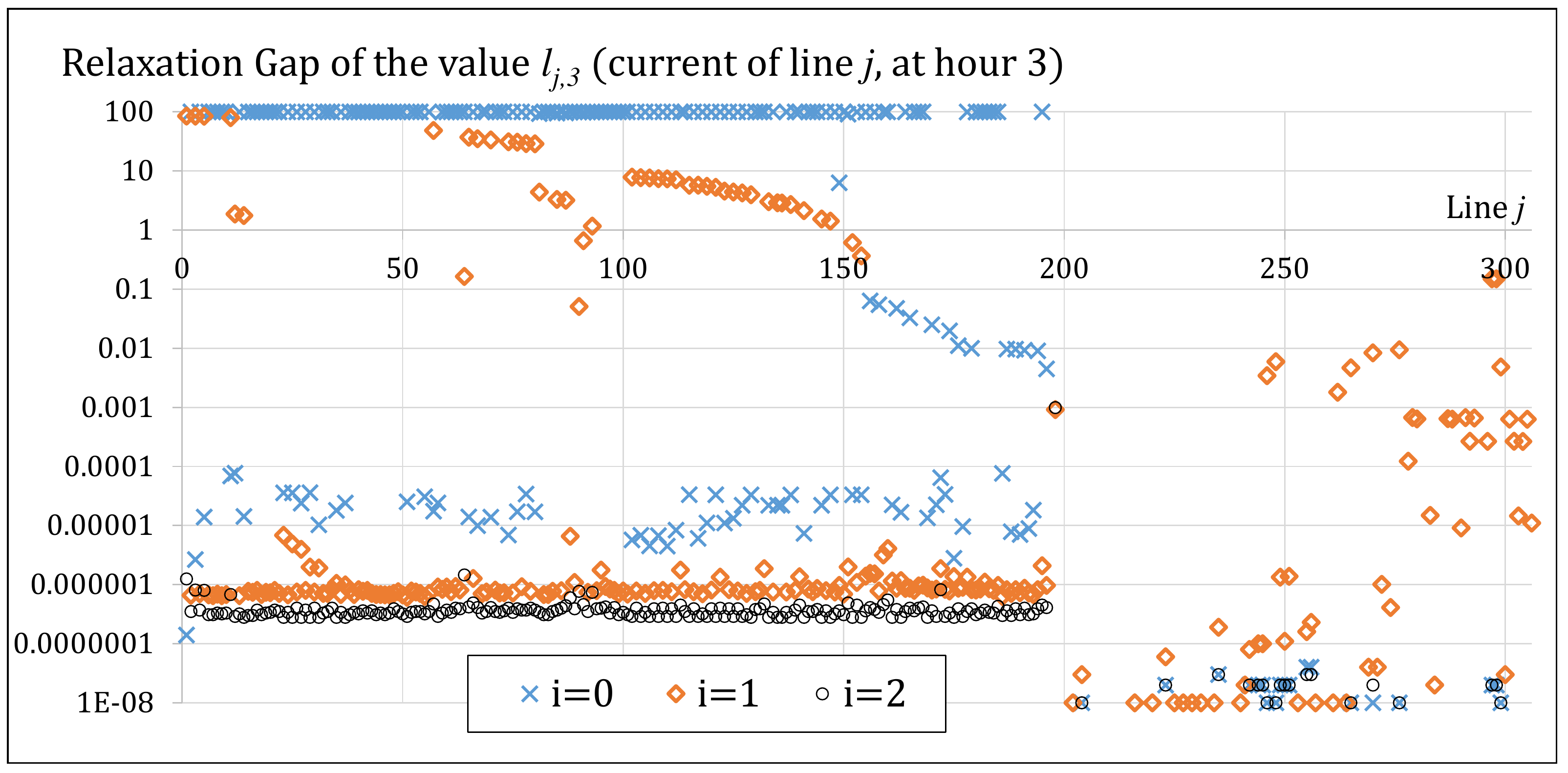}
\caption{Gap of calculated current values at hour 3; EV scenario, centralized solution, 3 iterations. Lines 197 to 306 (last 110) refer to transformers. Values below 1e$^{-08}$ are not shown.} 
\label{FigR6}
\end{figure}

Furthermore, We tested a case with a negative LMP 
(setting hour 3 LMP at -5 \$/MWh), 
which we know renders the SOCP relaxation not exact, 
and we applied the methodology proposed in Subsection \ref{NonExact}.
We selected $\bar \rho^0 = 0.005$,
and proceeded to multiply by a factor of $\frac{3}{2}$
that converged for $\tau = 10^{-4}$.
An interesting, albeit not surprising, observation is that 
initially only the currents in the distribution lines were inflated, 
while the service transformer currents exhibited no gap --- 
we anticipated this behavior 
since the transformer degradation in the objective function
penalizes excessive service transformer currents. 
This mitigates the perverse incentive occurring 
in non exact convex relaxation instances 
to artificially increase losses and profit by importing more of the negatively priced energy.
In Fig. \ref{FigR6}, 
we present the relaxation gap for the value of the current during hour 3.
The initial solution ($i= 0$) exhibits high gaps for the currents 
of actual distribution lines (first 196 values), 
whereas practically zero gaps for the transformers.
At iteration $i= 1$, the gaps for the lines are reduced,
but may increase for the transformers (not to a high value though).
At $i= 2$, 
all gaps drop to practically acceptable levels.

\subsubsection{{\color{black} Sensitivity Analysis of Reactive Power Opportunity Cost}}
{\color{black}
Our results so far assumed an opportunity cost for reactive power at the substation equal to 10\% the value of LMP.
We performed sensitivity analysis with respect to this parameter, 
and evaluated the proposed method under values that equal 2\%, 4\%, 6\%, and 8\%,
under all three scenarios of Subsection \ref{Results}, 
namely the \textbf{EV} scenario (662 EVs), 
the \textbf{PV} scenario (220 PVs), 
and the combined \textbf{EV+PV} scenario (662 EVs and 220 PVs).
The obtained results exhibited similar performance with the 10\% value considered in Subsection \ref{Results} 
--- see Fig. \ref{Fig9} which presents the convergence progress for the three scenarios and the different reactive power opportunity costs.

\begin{figure}[tb]
\centering
\includegraphics[width=3.45in]{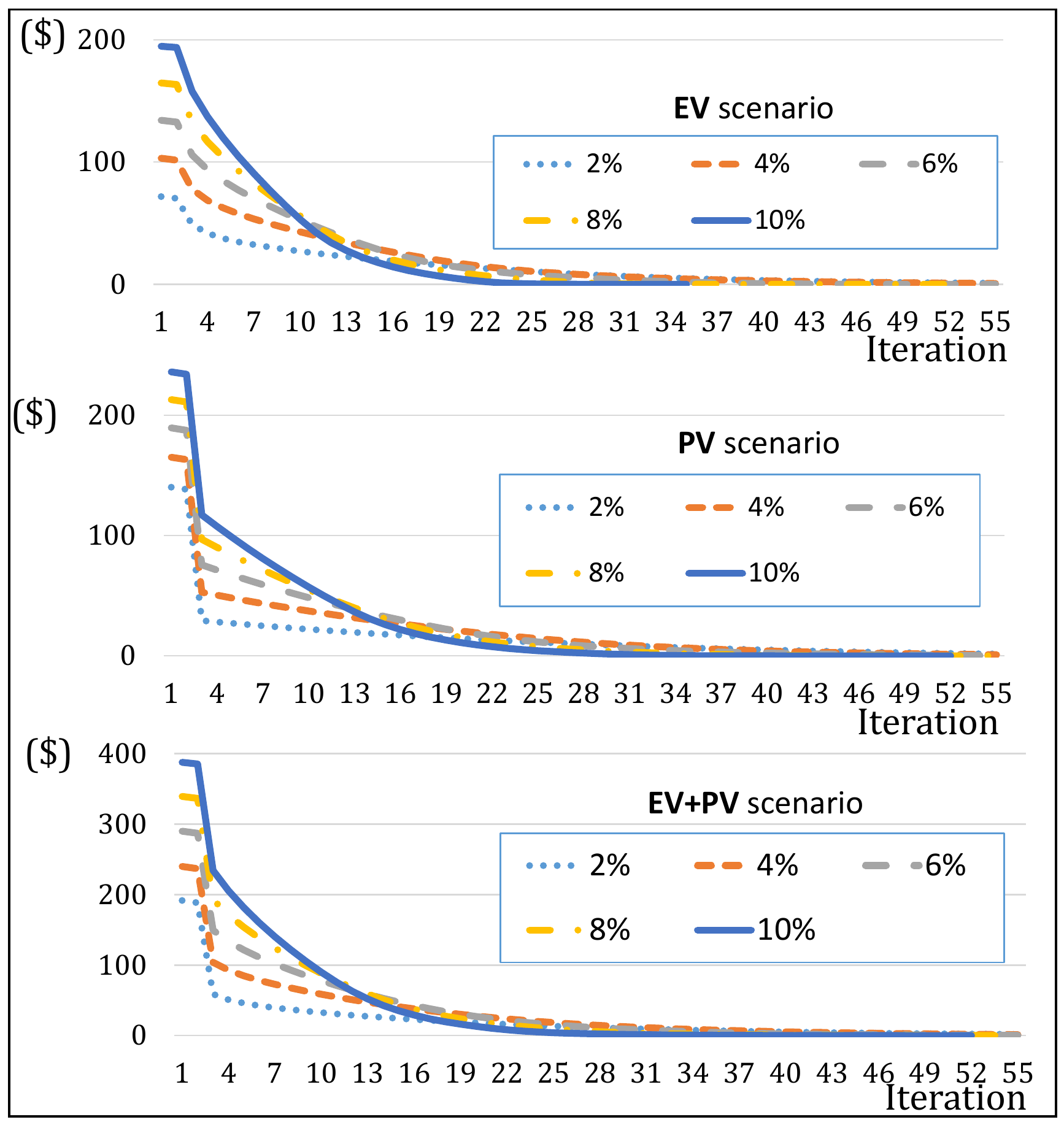}
\caption{{\color{black}Convergence progress of the 
decomposition method; sensitivity analysis \emph{w.r.t.} reactive power opportunity cost; \textbf{EV}, \textbf{PV}, and \textbf{EV+PV} scenarios. The x-axis shows the difference between the system cost at iteration $k$ and the system optimal cost (in \$).}} 
\label{Fig9}
\end{figure}

}

\subsubsection{{\color{black} Comparison with Time of Use tariffs}}
{\color{black}
Last but not least, we discuss the comparison of the proposed method 
with Time of Use tariffs and/or combination with marginal cost pricing.
The main problem of ToU tariffs is that they are not dynamic.
As such, ToU tariffs do not provide adequate signals for optimal DER scheduling.
Let us consider the example of the PV scenario, and the Q-DLMC trajectories at the service transformer locations shown in Fig \ref{Fig10} for the system optimal solution.
We observe that the Q-DLMCs differ at each service transformer may be above or below the red dashed line, i.e., higher or lower than the opportunity cost at the substation.
Hence, it is evident that a ToU tariff for reactive power, which would offer the same price at each location differentiated only with time,
would lose the locational signal,
and would induce the provision of a higher or lower (in any case sub-optimal) amount of reactive power.
Higher provision of reactive power may overload the transformer, increase marginal degradation and create spikes (as the ones observed in Fig. \ref{FigR7} (also discussed in our prior work \cite{IEEEtsg20}).
Lower provision of reactive power would result in sub-optimal system solutions.
Note that Q-DLMCs become also zero at certain locations and hours, however,
as we mentioned earlier the optimal quantities for reactive power are not zero; 
they are discovered by the proposed decomposition (leveraging the proximal terms).
We should also note that mixing price signals that depend on the marginal cost for real power but a fixed value (or ToU) for reactive power is not compatible with our proposed framework, because these marginal costs cannot be decoupled.
Our method captures this coupling, and produces price signals that reflect the marginal cost for both real and reactive power.
If we removed this coupling, the outcome would move away from the efficient solution (as also Fig. \ref{Fig10} suggests).
Hence, even if we provided price signals to DERs that reflected the optimal marginal costs for real power (suppose we could guess these costs) but not for reactive power,
then the outcome would be new (sub-optimal) DER quantities, 
which would drive also the real power marginal costs away from their optimal values.

\begin{figure}[tb]
\centering
\includegraphics[width=3.45in]{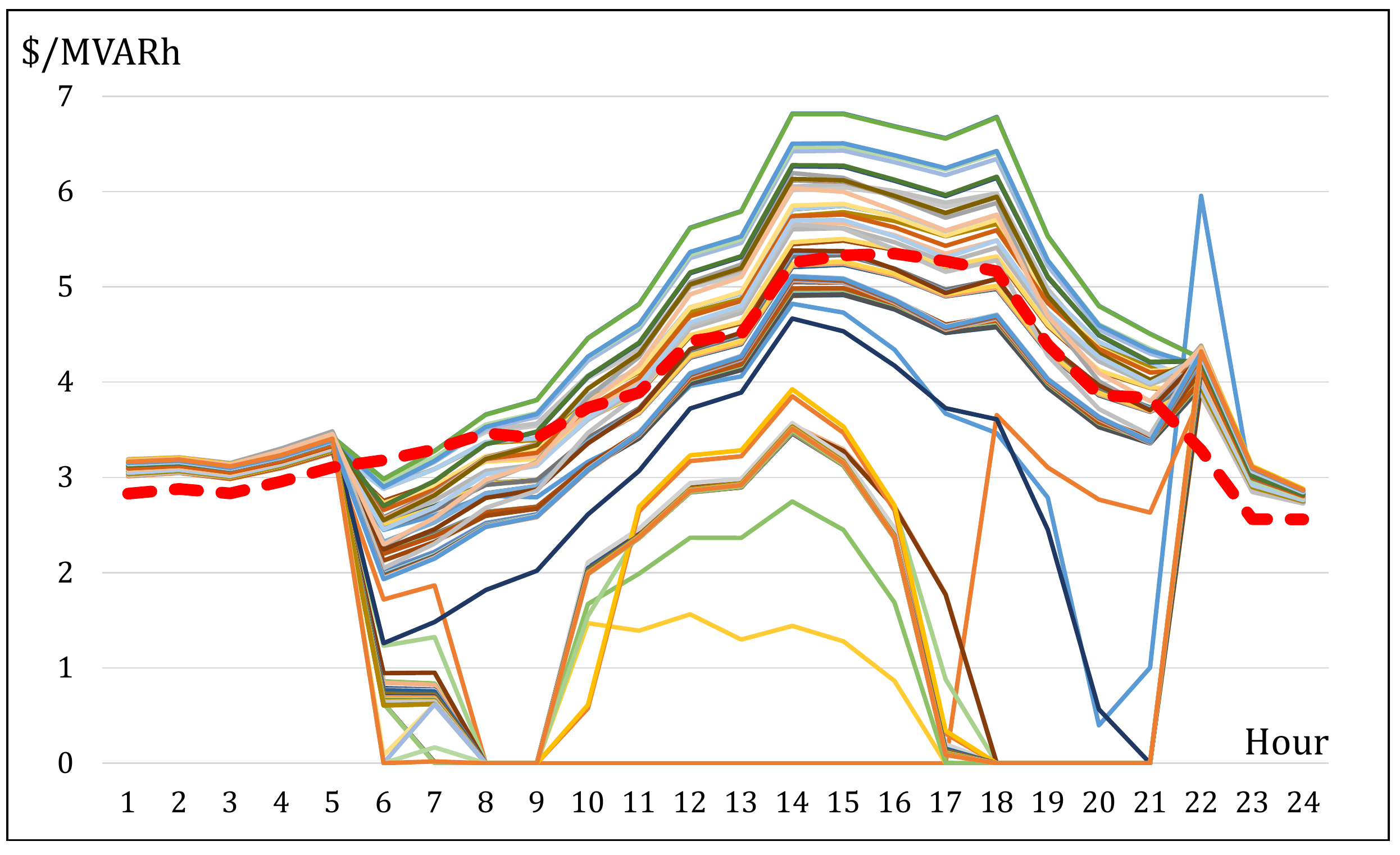}
\caption{{\color{black}Q-DLMC trajectories for each service transformer location; \textbf{PV} scenario; system optimal solution. The red dashed line shows the reactive power opportunity cost at the substation (10\% the value of LMP).}} 
\label{Fig10}
\end{figure}

}
{\color{black}
\section{Discussion} \label{Discuss}

The proposed 
decomposition offers a method for solving the Grid-DER coordination problem \eqref{SocialCost}--\eqref{DERCon}. 
This method has a natural economic interpretation. 
The DER {\color{black} optimization} problems are profit maximization problems (with an additional regularization term).
The tentative ``prices'' offered to the DERs at each iteration 
reflect the marginal cost calculated by the Grid Operator 
conditional upon the DER schedules determined and submitted at the previous iteration. 
When the DERs see the new tentative ``prices,'' they change, 
i.e., adapt their schedules to the new tentative prices; 
the adapted schedules are passed on to the Grid Operator,
who updates the spatiotemporal system marginal costs (``prices''), 
and the process repeats until convergence, 
i.e., until a fixed point of system spatiotemporal marginal costs and DER schedules. 
Most importantly, the marginal costs revealed upon convergence, 
i.e., when we obtain the optimal/fully adapted DER schedules, 
are indeed the marginal costs associated with the efficient operating point 
of the system that constitutes a mutually beneficial equilibrium. 
If these were the ``prices'' finally offered/cleared in a potential market setting, 
the equilibrium DER quantities would maximize DER profits 
and the system would minimize its costs. 
Hence, these revealed marginal costs (``prices'') 
can be used for the design of tariffs (and opt-in schemes as described in \cite{bor}) 
to move towards dynamic-pricing-based schemes 
and incentivize DERs to provide their full and truthful preferences and capabilities
to the platform which would enable the iterative derivation of a mutually beneficial equilibrium.
Dynamic spatiotemporal rates (as opposed to static flat or Time of Use rates) along with an opt-in scheme, 
will create what \cite{bor} calls a ``virtuous cycle,'' 
which will lead more and more customers to join the dynamic rates. 
In fact, inelastic customers who may opt for average cost prices 
will also benefit from progressively decreasing average costs
that will result from the decreasing marginal costs. 

We envision this adaptively dynamic scheme materializing through a platform, 
where DER owners submit their preferences 
(e.g., EVs submit their anticipated location and required charging deadlines for the next day etc.). 
One can think of this as a potential market design allowing ``complex bids/offers'' 
where DERs are not bound to uniform hourly price quantity bids/offers, 
but are free to provide the platform --- quite possibly anonymously --- 
their intertemporal constraints, preferences, requirements and costs, 
enabling their use as input to the proposed decomposition method. 
Note that DER information is provided once. 
In other words, DERs do not ``game'' by changing their costs/preferences at each iteration. 
Because they are too many and too small, and because they do not have network information, 
it is reasonable to assume that they will not collude or behave strategically. 
This, however, may not be the case when DER aggregators are given the authority to schedule DERs,
and especially when these aggregators have knowledge of the network costs structure and can anticipate the impact of DER schedules on System marginal costs. 
This is an issue discussed in prior work \cite{yac20}, 
where we have shown that network-information-aware DER aggregators, 
as for example the Distribution System Operator itself, 
may lead to inefficient market outcomes.

We believe we make a significant contribution 
with respect to the common perception on the anticipated fluctuation of marginal costs in distribution networks. 
Our analysis shows that high volatility of marginal costs is very significantly mitigated when DERs are coordinated, 
i.e., when they are able to adapt to spatiotemporal network marginal costs 
leading to optimal schedules constituting a mutually beneficial equilibrium. 
Extreme marginal cost instances are only associated with contrived 
and most importantly severely sub-optimal DER schedules that our adaptive framework avoids. 
We are of course cognizant of the fact that the DER schedules our process derives 
constitute a day ahead operational plan 
and that the actual conditions driving real time DER behavior will differ.
We can argue, however, 
that these differences that the network operator should anticipate and plan on, 
will be significantly reduced 
if the day ahead operational plan reflects the complex daily bids 
that each DER is allowed to submit, and, which, 
undoubtedly incorporate each DER's forecast of local weather (PV) and travel plans (EV). 
{\color{black} This expected reduction in uncertainty resulting from our ability to handle complex DER bids, 
reinforces the value of solving the apparently deterministic day ahead problem 
through our proposed 
{\color{black} decomposition method}. 

Indeed, as more and more stochastic resources are being integrated in the distribution grid,
the role of uncertainty is becoming increasingly significant.
There are several works that model AC OPF in a stochastic setting, 
e.g., through chance-constraints in centralized control \cite{krah17} 
and in data-driven approaches \cite{bb19} that emphasize on the operational aspect.
A chance-constrained AC OPF formulation is also employed in \cite{md20},
which accounts for small-scale generators as controllable DERs,
while treating all behind-the-meter DER as uncontrollable,
and includes pricing considerations with chance-constrained generation and voltage limits.
On the other hand, \cite{CaramanisEtAl2016} is the first work 
that refers to reserves in the context of distribution network marginal pricing.
Indeed, reserve products have been a long-standing solution,
which, however, has not been yet explored in distribution networks.
Arguably, reserve products, ideally with endogenously determined requirements,
can be part of a solution that will mitigate the impact of variability and uncertainty in distribution networks,
and ensure that the employed voltage/ampacity limits in the day-ahead scheduling problem guarantee a secure operation in real-time.

This paper emphasizes the role of the deterministic scheduling problem in distribution networks.
Although we acknowledge that stochasticity plays an important role,
we argue that rigorous analysis of the deterministic case is still lacking, 
and that this paper offers a method, 
which solves significant practical problems, that are not solved by other methods,
by successively improving feasible AC OPF solutions in a computationally tractable manner.
In addition,} as DER penetration increases, 
we argue that a scheduling problem for the next day is fundamental 
for mitigating the expected real time volatility and achieving cost-efficient solutions. 
The reason is that, with the increasing penetration of DERs, 
an increasing portion of the load (by load we loosely refer to the net load) 
will become flexible and schedulable.
Hence, being able to schedule this load will improve our load forecasts, 
an aspect that is often neglected.
The requisite presence of information platforms and active grids 
is expected to occur in the near future 
given the imminent role of smart meters and big data analytics. 
{\color{black} The proposed decomposition method can leverage these advances,
and offer a framework for an efficient day-ahead scheduling, 
that will contribute to a more predictable real-time operation.}

We have also taken several steps to show, in both this and prior works (e.g., \cite{IEEEtsg20}), 
that other popular, open-loop, pricing schemes, e.g., Time of Use rates, 
do not provide adequate incentives to lead DERs to an efficient day ahead schedule. 
In addition, we have shown that if we neglect the transformer degradation cost, 
and focus only on loss minimization, 
the resulting schedule is associated with high transformer overloads and high marginal costs. 
The proposed framework is based on the combined effect of real and reactive power. 
Reactive power is essential in distribution networks. 
Moreover it interacts with real power in determining the tradeoffs 
between a specific DER's real-reactive power output, 
and also in the determination of system costs. 
Indeed, real and reactive power cannot be decoupled 
without compromising the determination of the intertemporal marginal costs
and the optimal solution in a non trivial manner.

}

\section{Conclusions} \label{Conclusions}

We presented a 
decomposition {\color{black} method} 
to the hourly coordination of a large number of DERs 
connected at a distribution feeder 
transacting real and reactive power with the transmission network, 
while adapting to granular spatiotemporally varying DLMCs 
downstream of each service transformer. 
Most importantly, 
we have implemented the proposed algorithm on an actual distribution feeder
featuring 110 service transformers, 
and have shown that a large number of feeder connected EVs and PVs 
can adapt to DLMCs and reach an optimal coordination in a reasonably fast converging 
decomposition framework. 
Future work will focus on improving further the algorithm by, among others, 
tuning regularization term penalties adaptively in order to achieve rapidly converging implementations.

\section*{Acknowledgment}
The authors would like to thank Holyoke Gas and Electric for providing actual feeder data for the pilot study. 
This research paper benefited from the support of the FMJH
Program PGMO and from the EDF support to this program.

\ifCLASSOPTIONcaptionsoff 
  \newpage
\fi

\end{document}